\documentclass[11pt]{amsart}
\usepackage{enumerate}
\usepackage{amsfonts}
\usepackage{amssymb}
\usepackage{amsmath}
\usepackage{amsthm}
\usepackage{mathtools}
\usepackage{amstext}
\usepackage{tikz-cd}
\usepackage[left=1.2in,right=1.2in,top=1.4in, bottom=1.4in]{geometry}
\usepackage{graphicx}

\usepackage[section]{placeins}
\usepackage{caption}

\usepackage{xcolor}
\definecolor{cit}{HTML}{117733}
\definecolor{lin}{HTML}{0003d6}
\usepackage{hyperref}
\hypersetup{bookmarksnumbered,colorlinks=true,linkcolor=lin,citecolor=cit,urlcolor=black}

\DeclareMathOperator{\bl}{Bl}
\DeclareMathOperator{\Nef}{Nef}

\DeclareMathOperator{\Cox}{Cox}
\DeclareMathOperator{\Mov}{Mov}
\DeclareMathOperator{\Pic}{Pic}
\DeclareMathOperator{\NE}{NE}
\DeclareMathOperator{\NEb}{\overline{\NE}}

\DeclareMathOperator{\rk}{rank}
\DeclareMathOperator{\lcm}{lcm}
\DeclareRobustCommand{\stirling}{\genfrac\{\}{0pt}{}}

\newcommand{\fl}[1]{\lfloor#1\rfloor}
\newcommand{\cl}[1]{\lceil#1\rceil}

\newcommand{\ac}{\mathcal{A}}
\newcommand{\au}{\Phi}

\newcommand{\CC}{\mathbb{C}}
\newcommand{\Z}{\mathbb{Z}}
\newcommand{\R}{\mathbb{R}}
\newcommand{\PP}{\mathbb{P}}

\newcommand{\Q}{\mathbb{Q}}
\newcommand{\oO}{\mathcal{O}}

\newcommand\btr[1]{\ensuremath{\sum_{#1=0}^n (-1)^#1 {n\choose #1}}}
\newcommand\wtr[1]{\ensuremath{\sum_{#1=0}^5 (-1)^#1 {5\choose #1}}}

\newcommand{\e}{\equiv}
\newcommand{\m}{\pmod}

\newcommand{\vv}[1]{\lvert#1\rvert}
\newcommand{\vr}[1]{\langle#1\rangle}
\newcommand{\ra}{\to}

\newcommand{\mt}{\mapsto}

\newcommand{\ik}{^{-1}}
\newcommand{\ifff}{\Longleftrightarrow}

\newcommand{\dis}{\displaystyle}
\newcommand{\txt}[1]{\text{ #1 }}
\newcommand{\para}{\vskip 1em}

\theoremstyle{plain}
\newtheorem{theorem}{Theorem}[section]
\newtheorem{corollary}[theorem]{Corollary}
\newtheorem{prop}[theorem]{Proposition}
\newtheorem{lemma}[theorem]{Lemma}
\newtheorem{conj}[theorem]{Conjecture}

\theoremstyle{definition}
\newtheorem{defi}[theorem]{Definition}

\newtheorem{example}[theorem]{Example}
\newtheorem{app}[theorem]{Application}
\newtheorem{remark}[theorem]{Remark}
\newcommand{\pf}{{\em Proof}}

\begin{document}
\title[New Examples and non-examples of MDS when Blowing up Toric Surfaces]{New Examples and Non-examples of Mori Dream Spaces when Blowing up Toric Surfaces}
\author{Zhuang HE}
\address{Zhuang He. Northeastern University. 360 Huntington Ave, Boston, MA 02115}
\email{he.zhu@husky.neu.edu}

\begin{abstract}
We study the question of whether the blow-ups of toric surfaces of Picard number one at the identity point of the torus are Mori Dream Spaces. For some of these toric surfaces, the question whether the blow-up is a Mori Dream Space is equivalent to countably many planar interpolation problems. We state a conjecture which generalizes a theorem of Gonz\'{a}lez and Karu. We give new examples and non-examples of Mori Dream Spaces among these blow-ups. 
\end{abstract}
\maketitle

\section{Introduction}
\label{intro}
Mori Dream Spaces (MDS) were introduced by Hu and Keel in \cite{hu2000} as normal, $\Q$-factorial projective varieties such that 
\begin{enumerate}
	\item $\Pic(X)_\Q=N^1(X)_\Q$.
	\item $\Nef(X)$ is generated by finitely many semiample divisors.
	\item There exist finitely many small $\Q$-factorial modifications (SQMs) $g_i:X\dashrightarrow X_i$, $i=1,\cdots r$, such that every $X_i$ satisfies (1) and (2), and the cone of movable divisors
		$ \Mov(X)$ is the union $\bigcup_{i=1}^r g_i^*(\Nef(X_i)).$
	\end{enumerate}
	Mori Dream Spaces are important examples of varieties where Mori's program can be run for every divisor \cite[Prop. 1.11]{hu2000}.

Mori Dream Spaces are related to Cox rings. Assume $X$ is a projective variety over the complex numbers $\CC$, with finitely generated Picard group. Choose line bundles $L_1,\cdots, L_s$ which span $\Pic(X)_\Q$. A Cox ring of $X$ is defined as the direct sum
\[\Cox(X):=\bigoplus_{(m_1,\cdots,m_s)\in \Z^s} H^0 (X, m_1 L_1+\cdots +m_s L_s).\]
The finite generation of $\Cox(X)$ does not depend on the choice of $L_1,\cdots,L_s$ \cite{hu2000}. Furthermore, it is shown in \cite[Prop. 2.9]{hu2000} that a  $\Q$-factorial projective variety $X$ with $\Pic(X)_\Q=N^1(X)_\Q$ is a MDS if and only if $\Cox(X)$ is a finitely generated $\CC$-algebra.

One basic source of examples of Mori Dream Spaces is from toric varieties. However, the blow-ups of MDS can fail to be MDS (as we recall below). For instance, blow-ups of toric varieties, in particular weighted projective planes, can fail to be MDS. A weighted projective plane $\PP(a,b,c)$ is the quotient of $\CC^{3}-\{0\}$ by the following $\CC^*$-action:
\begin{align*}
	\CC^*\times (\CC^{3}-\{0\})\ra \CC^{3}-\{0\}\\
	(t,(x,y,z))\mt (t^a x, t^b y, t^c z),
\end{align*}
where $a,b$ and $c$ are positive integers. Note that $\PP(a,b,c)$ is a toric projective surface of Picard number one.

We denote $X=\PP(a,b,c)$. Let $X'=\bl_e \PP(a,b,c)$ be the blow-up of $X$ at the identity point $e$ of the open torus.

\para{\bf Question:} For which $a,b,c$ is the blow-up $X'=\bl_e \PP(a,b,c)$ a MDS?
\para

Historically, Cutkosky gave many sufficient conditions for $\Cox(X')$ to be finitely generated, and equivalently, $X'$ to be a MDS. For example, if $-K_{X'}$ is big, then $X'$ is a MDS \cite{cutkosky1991}. In particular, if $a+b+c>\sqrt{abc}$ (for example, when one of $a,b,c$ is at most $4$), then  $-K_{X'}$ is big, and $X'$ is a MDS.  
Based on the work of Cutkosky, Srinivasan \cite{HS} attained several numerical conditions for $X'$ to be a MDS, including that if one of $a,b,c$ is $6$, then $X'$ is a MDS.
Recently, Hausen, Keicher and Laface gave an algorithm which provides more examples of MDS \cite{hkl}.

On the other hand, by 2013 the only known examples of blow-ups at a general point of toric varieties failing to be a MDS were given by Goto, Nishida and Watanabe \cite{goto1994} (1994). For example, they showed that $X'=\bl_e \PP(a,b,c)$ is not MDS when $(a,b,c)=(7N-3,8N-3, (5N-2)N)$ for $N\geq 4$, and $3\nmid N$. 
In 2013, Castravet and Tevelev \cite{castravet2015} proved, using the \cite{goto1994} results, that blow-ups of Losev-Manin moduli spaces at a general point are not MDS in sufficient large dimensions. Using this, they proved that the moduli space of curves $\overline{M}_{0,n}$ is not a Mori Dream Space for $n>133$.
In 2014, Gonz\'{a}lez and Karu \cite{GK} provided more examples of triples $(a,b,c)$ such that $X'$ is not a MDS, and lowered the bound to $n>12$, which was further improved to $n>9$ by Hausen, Keicher and Laface in \cite{hkl}.

In this paper, we further develop and generalize Gonz\'{a}lez and Karu's idea, and show that for some toric surfaces $X$ of Picard number one, the blow-up $X'$ at the identity point of the torus is not a MDS if and only if a family of countable many planar interpolation problems in $\PP^2$ all have solutions. (Proposition \ref{intpblm}, Corollaries \ref{int} and \ref{int2}). As an application, we provide new examples (Theorem \ref{main01}) and new non-examples (Theorems \ref{main03}, \ref{main04}) of MDS, from blow-ups of toric surfaces of Picard number one, in particular blow-ups of weighted projective planes (Example \ref{exdd}). Our method is different from Hausen, Keicher and Laface's. The results above can be combined into a numerical criterion (Corollary \ref{mainwidth}). We make a conjecture (Conjecture \ref{conj00}) generalizing the non-examples of \cite{GK}.

\para
{\bf Acknowledgement} The author is very grateful to Ana-Maria Castravet for introducing this topic and the guidance throughout the research. The author would like to thank Kalle Karu and Xiaolei Zhao for valuable correspondences and suggestions. The author thanks Jos\'{e} Gonz\'{a}lez for help with completing the proof of Proposition 5.1.

\tableofcontents

\section{Main Results}
\label{mainresult}

\para{\bf Notations and Settings} We work over $\CC$.
Let $N= \Z^2$ be the plane lattice, and $N_\R:=N\otimes_\Z \R \cong \R^2$.

A triangle $\Delta\subseteq N_\R$ with rational slopes $s_1, s_2$ and $s_3$ determines a polarized toric variety $(X_\Delta, H)$, in the way that the fan of $X=X_\Delta$ is the normal fan of $\Delta$, and $H$ is the $\Q$-Cartier torus-invariant divisor determined by $\Delta$. 
Notice that two triangles $\Delta$ and $\Delta'$ determine the same $X$ if they have the same slopes $s_1,s_2,s_3$.
Hence, without specific mention, we assume every triangle $\Delta$ is at the following position:
one vertex is at $(0,0)$, and the opposite side of $(0,0)$ passes through $(0,y_0)$ for $y_0>0$. 

Define $\Delta_0$ as the triangle with $y_0=1$.

A triangle in $N_\R$ is a {\em lattice triangle} if all its vertices are in $N$.
We say a lattice triangle $\Delta$ is {\em good} if in addition the $y$-intercept $y_0$ is an integer. 
A good lattice triangle $\Delta$ is an integer multiple of the $\Delta_0$ with the same slopes, and there exists a smallest good lattice triangle with the given slopes, which we denote as $\Delta_1$.
Let $m$ be the integer such that $\Delta_1=m\Delta_0$.

A {\em column} of a lattice triangle $\Delta$ is the set $\Delta\cap \{(x,y)\in \Z^2\mid x=x_0\}$ for some integer $x_0$, such that $x_0$ does not equal the $x$-coordinates of the leftmost and rightmost vertices of $\Delta$.

Given $s_1<s_2<s_3\in\Q$, the {\em width} of $s_1,s_2,s_3$ is
\[w=\frac{1}{s_2-s_1}+\frac{1}{s_3-s_2},\]
which equals the width of $\Delta_0$. Let $W:=kmw$, which is the width of $k\Delta_1$.
\para
{\bf Assumption: $w<1$. }
\para

Finally, $X'=\bl_e X_{\Delta}$ is the blow-up at the torus identity point $e$ of $X$. 

\begin{prop} $($see \cite[Prop. 8.6]{castravet16}$)$
	Given $s_1<s_2<s_3$ rational numbers with $w<1$, Let $X$ the toric variety defined by $s_1,s_2,s_3$. Let $X'$ be the blow-up of $X$ at the torus identity point $e$. Then the blow-up $X'$ is not a MDS if and only if for every sufficiently divisible integer $k>0$, there exists a curve $Y$ in $\PP^2$, of degree up to $W-1=kmw-1$, and a vertex $p\neq (0,0)$ of $\Delta_1$,  such that $Y$ passes through all the points $(i,j)\in k\Delta_1\cap \Z^2$ but does not pass through $kp$.
	\label{intpblm}
\end{prop}

\pf. See Section \ref{cox}.\qed

\para
We consider the interpolation problem proposed in Proposition \ref{intpblm}. For any $\Delta=k\Delta_1$, the column at $\{x=0\}$ contains $km+1$ points. Suppose such a curve $Y$ in degree $\leq kmw-1=W-1<km+1$ exists, then $Y$ passes through all those $km+1$ points. By B\'{e}zout's Theorem, $Y$ must contain the line $\{x=0\}$ as a component. Hence the existence of $Y$ is equivalent to the existence of $Y_1$ of degree $\leq kmw-2$, passing through all integer points except the column at $\{x=0\}$, and $kq$, but not passing through $kp$. 

Indeed this argument by B\'{e}zout's Theorem can be run for the rest of columns in $k\Delta_1$ until step $n$, where $n$ is the integer such that every remaining columns in $k\Delta_1$ after step $n$ contain no more than $kmw-1-n=W-1-n$ points. Notice that when we stop, the existence of $Y$ is equivalent to the existence of a curve $Y_n$ of degree $\leq W-1-n$, passing through all columns remaining, and $kq$, but not passing through $kp$.

\begin{defi}\label{redd}
	The {\em reduced degree} of a good lattice triangle $\Delta$ of width $W$ equals $W-1-n$, where $n$ is the maximal number of steps we can run as above.
	Equivalently, $d$ equals the number of remaining columns (the left and right vertices excluded by definition) in $\Delta$ after we deleted all columns through the above reduction process.
\end{defi}

We now introduce the reduced degree $d$ and minimal degree $d'$ of a triple of slopes $s_1,s_2,s_3$.

\begin{defi}	
	\label{redd2}
	Given rational numbers  $s_1<s_2<s_3$, with $w<1$, let $\Delta_1$ be the smallest good lattice triangle with slopes $\{s_i\}$. The {\em reduced degree} of $s_1,s_2,s_3$ is the largest nonnegative integer $d$, such that there are exactly $d$ columns in $\Delta_1$ containing $\leq d$ points.
\end{defi}
\begin{defi}	
	\label{mind}
	Given rational numbers  $s_1<s_2<s_3$, with $w<1$, let $\Delta_1$ be the smallest good lattice triangle with slopes $\{s_i\}$. The {\em minimal degree} of $s_1,s_2,s_3$ is the smallest positive integer $d'$, such that there are exactly $d'$ columns in $\Delta_1$ containing $\leq d'$ points. When no such $d'>0$ exists, we define the minimal degree to be zero.
\end{defi}

\begin{theorem}\label{unity}
	Consider any rational numbers $s_1<s_2<s_3$ such that the width $w<1$. 
	\begin{enumerate}
		\item The reduced degrees of the good lattice triangles $\Delta$ with slopes $s_1,s_2,s_3$ all equal the reduced degree of $s_1,s_2,s_3$ (Definition \ref{redd2}). In particular, the reduced degree of $s_1,s_2,s_3$ equals the reduced degree of $\Delta_1$.
		\item The reduced degree $d$ satisfies the following inequality:
			\[d<\frac{w}{1-w}.\]
	\end{enumerate}
\end{theorem}

\begin{figure}[h]
	\includegraphics[width=4cm]{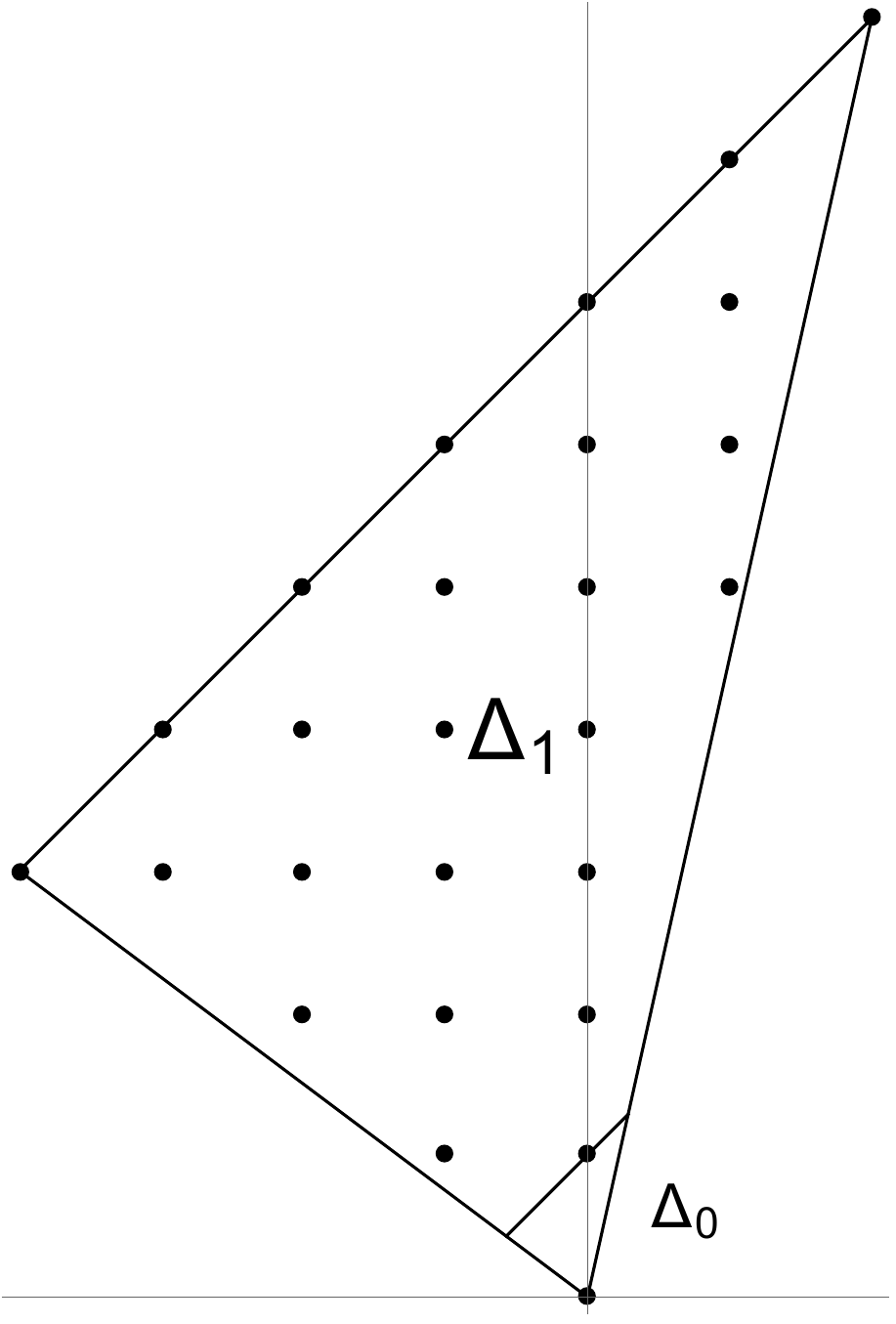}
	\caption{The triangle $\Delta_1$ with slopes $(-3/4,1,9/2)$ has reduced degree $0$}
	\label{fig:uu}
\end{figure}

For example, the triangle $\Delta_1=7\Delta_0$ with slopes $(-3/4,1,9/2)$ is shown in the Figure \ref{fig:uu}. There are $5$ columns in $\Delta_1$, at $x=-3,-2,-1,0,1$. The number of points in each column is $2,4,6,8,4$. The numbers of columns with $\leq i$ points are given below for $0\leq i\leq 8$.
As a result, the reduced degree $d=0$. Further, for any integer multiple $k\Delta_1$, the reduced degrees will not increase, and all equal to $0$. Finally, the minimal degree $d'$ equals $0$.

\begin{table}[h]
	\begin{tabular}{|c|ccccccccc|}\hline
		$i$ &0&1&2&3&4&5&6&7&8\\\hline
	number of columns with $\leq i$ points& 0 & 0 & 1 & 1& 3 & 3 &4 &4 &5\\ \hline
	\end{tabular}
\end{table}

More interesting examples are the triangles in Figure \ref{fig:non} which have minimal degrees $d'=3,5,7$ respectively. (Recall that we do not include the leftmost and rightmost vertices as columns)

\begin{theorem}
	Consider rational numbers $s_1<s_2<s_3$ with width $w<1$. Let $d$ be the reduced degree and $d'$ be the minimal degree. Let $X$ the toric variety defined by $s_1,s_2,s_3$. Let $X'$ be the blow-up of $X$ at the torus identity point $e$.  
	\label{main01}
\begin{enumerate}
		\item If $d=0$, then the blow-up $X'$ is a Mori Dream Space {\em (}MDS\text{\em )}.
		\item If $d'=1$, then $s_2\not\in \Z$, and the blow-up $X'$ is not a MDS. 
		\item In particular,  if $d=1$, then $X'$ is not a MDS.
	\end{enumerate}	
\end{theorem}

For the proofs of Theorems \ref{unity} and \ref{main01}, see Section \ref{rd}.

Since $X'$ is always a MDS when $d=0$, it is worthwhile to develop a criterion for having $d=0$. In Section \ref{4} we prove the following combinatorial criterion:
\begin{prop}
	Consider rational numbers $s_1<s_2<s_3$ with width $w<1$. Let $X,X'$ be defined as in Theorem \ref{main01}. Let
	\[ l_k=\fl{k s_2}-\cl{k s_1}+1,\quad r_k=\fl{k s_3}-\cl{k s_2}+1,\quad \txt{for} k=1,2,\cdots\]
	Let $\gamma$ be the smallest positive integer such that
			$\gamma s_2^2, \gamma s_3$ and $\gamma s_2 s_3$ are all integers. 
			Then the reduced degree $d$ equals $0$ if all the following three conditions hold:
	\begin{enumerate}
		\item $l_1\neq 1$ and $r_1\neq 1$;
		\item If $\gamma>1$, then 
			\[s_2-s_1\geq \max_{1\leq t\leq \gamma-1} \frac{r_t+\{(r_t-t)s_2\}}{r_t-t};\]
		\item
			\[s_2-s_1\geq \frac{\gamma (s_3-s_2)+1+\{s_2\}}{\gamma(s_3-s_2-1)+1}\]
	\end{enumerate}
	\label{main02}
	where $\{x\}=x-\fl{x}$ is the fractional part of $x$. In particular, $X'$ is a MDS when all the conditions above hold.

	Conversely, if in addition we assume 
	\begin{enumerate}\setcounter{enumi}{3}
		\item For some $t$ where $\displaystyle \frac{r_t+\{(r_t-t)s_2\}}{r_t-t}$ achieves its maximum among $1\leq t\leq \gamma-1$, $t$ satisfies
\[\gamma(s_3-s_2-1)(r_t-1+\{(r_t-t)s_2\})-(r_t-t)(\gamma (s_3-s_2)+1)\leq 0,\]
	\end{enumerate}

Then $d=0$ implies (1), (2), (3) all hold.
\end{prop}

\begin{corollary}\label{zs2}
	Consider rational numbers $s_1<s_2<s_3$ with width $w<1$. If $s_2\in\Z$, then the reduced degree $d$ equals zero, and $X'$ is a MDS.
\end{corollary}
\pf. It suffices to assume $s_2=0$ (See Remark \ref{shear}). Then $l_k=-\cl{ks_1}+1$ and $ r_k=\fl{k s_3}+1$. Since $w=1/s_3-1/s_1<1$, we have $s_3>1$ and $s_1\leq -1$. So $l_1\geq 2$ and $r_1\geq 2$. This shows (1) of Proposition \ref{main02}. 

For (2), we find $\{(r_t-t)s_2\}=0$. We claim that for every $t\geq 1$, $r_t/(r_t-t)\leq s_3/(s_3-1)$. Indeed this simplifies to $r_t\geq s_3 t$, which holds since $r_t=\fl{s_3 t}+1\geq s_3 t$. Now $w<1$ implies that $-s_1\geq   s_3/(s_3-1)$. Therefore (2) holds.

Finally (3) follows from the inequality
\[-s_1\geq\frac{s_3}{s_3-1}> \frac{\gamma s_3+1}{\gamma(s_3-1)+1}.\]
Now Proposition \ref{main02} shows that $d=0$, so $X'$ is a MDS.\qed

\begin{corollary}
	For any $s_2, s_3\in \Q$, with $\fl{s_3}-\cl{s_2}\geq 1$, there exists a rational number $s_0$, depending on $s_2$ and $s_3$ only, such that if $s_1<s_0$, then the blow-up $X'$ is a MDS.
\end{corollary}

\pf. Firstly, $\fl{s_3}-\cl{s_2}\geq 1$ is equivalent to $r_1\geq 2$. If in addition $w<1$ and (1)-(3) of Proposition \ref{main02} all hold, then $X'$ is a MDS. The condition $w<1$ gives $s_1<(s_2 s_3-s_2^2-s_3)/(s_3-s_2-1)$. Assuming $w<1$, $l_1$ cannot be zero. Then $l_1\geq 2$ is equivalent to $s_1\leq \fl{s_2}-1$. Now $s_0$ can be taken as the minimum of the $4$ upper bounds above from $w<1$, $l_1\geq 2$, and (2)(3) of Proposition \ref{main02}. \qed

For example, consider $s_2=1/2$, and $s_3=3$ (see Example \ref{a5}). Here $w<1$ if and only if $s_1<-7/6$. $l_1\geq 2$ if and only if  $s_1\leq -1$. 
In Proposition \ref{main02}, we have $\gamma=4$. Then (2) says $s_1\leq -6/5$, and (3) says $s_1\leq -8/7$. Therefore we can take $s_0=-6/5$, so that when $s_1<-6/5$, $s_2=1/2$ and $s_3=3$, the blow-up $X'$ is a MDS.

For the case when the minimal degree $d'$ satisfies $d'\geq 2$ (hence, $d\geq 2$), we have the following conjecture:

\begin{conj}
		Consider rational numbers $s_1<s_2<s_3$ with width $w<1$. Let $d'$ be the minimal degree. Let $X$ the toric variety defined by $s_1,s_2,s_3$, and $X'$ be the blow-up of $X$ at the torus identity point $e$.  
If $d'\geq 2$, and $d'\cdot s_2\not\in\Z$, then the blow-up $X'$ is not a MDS.
		\label{conj00}
\end{conj}

\begin{remark}\label{gknonex}
	This conjecture, together with Theorem \ref{main01} (2), generalizes Gonz\'{a}lez and Karu's non-examples. Recall that in \cite{GK} Gonz\'{a}lez and Karu showed that if the lattice triangle $\Delta$ satisfy that
	\begin{enumerate}
	\item The first column from left has $n$ points;
	\item The $i$-th column from the right have $i+1$ points, for $i=1,2,\cdots,n-1$.
\end{enumerate}
Then the blow-up $X'$ is not a MDS if $ns_2\not\in\Z$.
Indeed, here the triangle $\Delta$ has minimal degree $d'=n$. 
Also note that Theorem \ref{main01} (2) is exactly the case of $n=1$.
\end{remark}

\begin{remark} The main observation is that we can classify all possible triangles $\Delta_1$ of a given minimal degree $d'$, by the numbers of lattice points on the columns with $<d'$ points.

Suppose $d'\geq 2$.
Since $d'\neq 1$, Lemma \ref{strictI} implies that the numbers of points on each columns are strictly increasing. Now there are exactly $d'$ columns in $\Delta_1$ with $\leq d'$ points, so those $d'$ columns must have $2,3,\cdots,d'-1,d',d'$ points respectively. 
Hence $\Delta_1$ determines a partition $\{2,\cdots, d'-1\}=S\sqcup T$, such that the number of lattice point on the columns starting from the left (right) vertex are given by $S$ ($T$ respectively). In particular, Gonz\'{a}lez and Karu's non-examples are given by $S=\emptyset$ and $T=\{2,\cdots,d'-1\}$.

As a result, it is helpful to classify all possible triangles $\Delta_1$ of a given minimal degree $d'$, by the numbers of lattice points on the columns with $< d'$ points, which are given by a partition  $\{2,\cdots, d'-1\}=S\sqcup T$, up to a horizontal reflection about the $y$-axis.

\end{remark}

With the assistance of computer programs (in {\em Mathematica 10} \cite{mma104}), we have:
\begin{theorem} \label{main03} Conjecture \ref{conj00} holds for $d'\leq 9$. That is,  If $2\leq d'\leq 9$, and $d'\cdot s_2\not\in\Z$, then the blow-up $X'$ is not a MDS.
	
\end{theorem}
\begin{theorem}\label{main04}
	Consider rational numbers $s_1<s_2<s_3$ with width $w<1$. Let $d'$ be the minimal degree.
	If $2\leq d'\leq 9$ and $d'\cdot s_2\not\in \Z$, then either we are in the case of Gonz\'{a}lez and Karu's non-examples \cite{GK} (Remark \ref{gknonex}), or $d'=5,7$ or $9$, and up to adding a same integer $t$ to all the three slopes, and a reflection about the $y$-axis, their slopes satisfy one of the following:
	\[\begin{cases}
		d'=5;\\
		-2-\frac{1}{2}<s_1\leq -2;\\
		\frac{1}{3}<s_2<\frac{1}{2}; \\
		2\leq s_3<2+\frac{1}{3}.
	\end{cases},
	\begin{cases}
		d'=7;\\
		-2-\frac{1}{3}<s_1\leq -2;\\
		\frac{1}{4}<s_2<\frac{1}{3}; \\
		2\leq s_3<2+\frac{1}{4}.
	\end{cases},
	\begin{cases}
		d'=9;\\
		-2-\frac{1}{4}<s_1\leq -2;\\
		\frac{1}{5}<s_2<\frac{1}{4};\\
		2\leq s_3<2+\frac{1}{5}.
	\end{cases}\]
respectively. Conversely, any combination of slopes which satisfies one of the system of inequalities above as well as $w<1$ and $d'\cdot s_2\not\in\Z$ (equivalently, $d'\cdot s_2\neq 2$) determines a blow-up $X'$ which is not a MDS. Figure \ref{fig:non} shows the relative positions of the lattice points on the columns with at most $d'$ lattice points in each lattice triangle.
\begin{figure}[h]
	\centering
	\includegraphics[width=\textwidth]{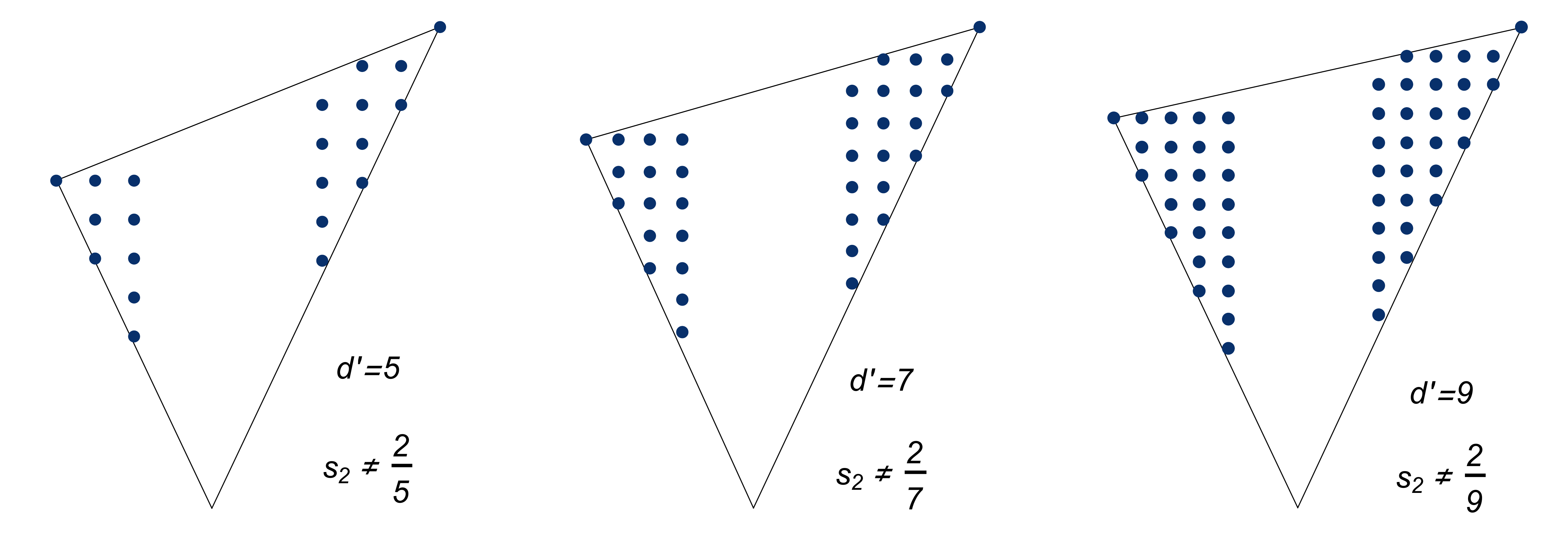}
	\caption{New non-examples of MDS given by the position of lattice points at corners}
	\label{fig:non}
\end{figure}

\end{theorem}

We will prove Theorem \ref{main03} and \ref{main04} in Section \ref{dayu2} and Section \ref{03pf}.

The inequalities of Theorem \ref{main04} come from the following general result.
\begin{lemma}\label{shape} Let $\Delta_1$ be the smallest good lattice triangle with slopes $s_1<s_2<s_3$. Let $n\in \Z_{>0}$. 
	Suppose in $\Delta_1$ the first $n$ columns from the left have $3,5,\cdots,2n-1,2n+1$ points, and the first $n+1$ columns from the right have $2,4,\cdots,2n,2n+1$ lattice points. Then up to adding a same integer $t$ to all the three slopes, and a reflection about the $y$-axis, the slopes satisfy the inequalities
	\[\begin{cases}
		\dis-2-\frac{1}{n}<s_1\leq -2;\\
		\dis\frac{1}{n+1}<s_2<\frac{1}{n};\\
		\dis 2\leq s_3<2+\frac{1}{n+1}.
	\end{cases}\]
\end{lemma}

\pf. See Section \ref{03pf}.\qed

Combining Theorem \ref{unity} (2), \ref{main01}, \ref{main03} and \ref{main04}, we have the following corollary.
\begin{corollary}
	Give rational numbers $s_1<s_2<s_3$ with width $w<1$. Let $d,d',X,X'$ be defined as in Theorem \ref{main01}. Then
	\begin{enumerate}
		\item If $w\leq 1/2$, then $d=0$, and $X'$ is a MDS.
		\item If $w\leq 10/11$ and $d'\cdot s_2\not\in\Z$,
then $1\leq d'\leq d\leq 9$, and $X'$ is not a MDS.
\end{enumerate}
	\label{mainwidth}
\end{corollary}

\begin{remark}\label{shear}
	 Indeed, adding the same integer $l$ to all slopes is equivalent to a shear transformation $(x,y)\mapsto (x,y+lx)$ of the triangles, which induces one-to-one correspondences between the solutions of the corresponding interpolation problems, and
 further induces isomorphisms between the corresponding toric varieties.
\end{remark}

\para
\begin{app}\label{appwpp}
	We apply Theorem \ref{main01} and Proposition \ref{main02} to weighted projective planes $\PP(a,b,c)$ in Section \ref{wpp}.

	By \cite[Lem. 10]{cutkosky1991}, for any triple $(a,b,c)$, there exist another triple $(a',b',c')$ of positive integers, such that $a',b',c'$ are pairwise coprime, and $\bl_e \PP(a,b,c)$ is isomorphic with $\bl_e \PP(a',b',c')$. Therefore we can assume that $a,b,c$ are pairwise coprime.
For every pairwise coprime triple $(a,b,c)$, there exists at most one triple of integers $(e,f,-g)$ such that $e,f,g>0$, $ae+bf=cg$, and $\gcd(e,f,g)=1$, even after permuting $a,b$ and $c$ (See \cite{GK}). We say such $(e,f,-g)$ a {\em relation} of the triple $(a,b,c)$, and always rearrange the triple $(a,b,c)$ in the order such that $ae+bf=cg$.

In Example \ref{exdd} we provide two tables of the complete lists of triples $(a,b,c)$ with $a\leq 15$ and $w<1$ such that the minimal degree $d'=1$ or $d'\geq 2$. In other words, any pairwise coprime triples with $a\leq 15$ that do not appear in any of the tables have $d=d'=0$, and give blow-ups which are MDS.
\end{app}

\begin{remark}
	Compared with the examples given in \cite{hkl}, we find that when $a,b,c\leq 30$, our examples of MDS form a proper subset of their examples. Specifically, when $a,b,c\leq 30$ and are pairwise coprime, there are four different cases:
	
	(1) There exists no relation such that $w<1$, even after permuting $a,b$ and $c$.

	(2) There exists a relation with $w<1$. The reduced degree $d\geq 1$ and it belongs to Gonz\'{a}lez and Karu's non-examples.

	(3) There exists a relation with $w<1$.
The reduced degree $d\geq 1$, and it does not belong to Gonz\'{a}lez and Karu's non-examples. There are only four of them $(8,13,15)$, $(8,13,25)$, $(15,19,29)$ and $(15,26,29)$.
	\begin{table}[h]
		\centering
		\begin{tabular}{l|cccc}\hline
			$(a,b,c)$ &$(8,13,15)$ & $(8,13,25)$ & $(15,19,29)$ & $(15,26,29)$ \\\hline
			minimal degree&3&  2 & 3& 4\\\hline
			reduced degree&3&  2 & 3& 4\\\hline
			$s_2$ & 2/3 &  1/2 & 2/3 & 3/4 \\\hline
		\end{tabular}
	\end{table}
	For those four triples, $d=d'$ and $d'\cdot s_2\in\Z$ (See the table below).

	(4) There exists a relation with $w<1$ and $d=0$, hence giving MDS blow-ups.

	It is worth pointing out that the examples of MDS by \cite{hkl} not only include all triples in cases (3) and (4), but also some in (1), for example, $\PP(7,18,19)$.
\end{remark}

On the other side, Theorem \ref{main04} provides the following new non-examples, all of minimal degree $d'=5$, including:
	\begin{enumerate}
	\item $a=11,f=1,g=7$, such that $49/11<b/c< 41/9$. The first examples are $(11,58,13), (11,140,31),(11,157,35),\cdots$.
	\item $a=13,f=1,g=8$, such that $64/13<b/c< 46/9$. The first examples are $(13,84,17), (13,149,30),(13,157,31),\cdots$.
	\item $a=14,f=1,g=9$, such that $81/14<b/c< 53/9$. The first examples are $(14,181,31), (14,309,53),(14,331,57),\cdots$.
	\item $\cdots$
\end{enumerate}

\section{The Interpolation Problem}
\label{cox}

We recall the following set-up from \cite{GK}. Fix rational numbers  $s_1<s_2<s_3$. Consider the unique triangle $\Delta_0$ given by slopes $(s_1,s_2,s_3)$, such that one vertex of $\Delta_0$ at $(0,0)$ and the opposite side of $(0,0)$ passes through $(0,1)$.
Let $(X,H)$ be the polarized toric surface determined by $\Delta_0$. That is, the normal fan of $\Delta_0$ is the fan of $X$, and $H$ is the $\Q$-Cartier divisor corresponding to the polytope $\Delta_0$. Then it follows that $H^2=w$, which equals to twice the area of $\Delta_0$. 

Now $\Pic(X)\otimes \Q$ is generated by $H$. Let $\pi: X'\ra X$ be the blow-up of $X$ at the torus identity point $e$. Let $E$ be the exceptional divisor of $X'$. Then $\Pic(X')\otimes \Q$ is two-dimensional, generated by the classes of $H'=\pi^* H$ and $E$. :

The section $1-y$ in $\oO_X(H)$ defined by the lattice points $(0,0)$ and $(0,1)$ (both on sides of $\Delta_0)$ gives us an effective divisor in the class $[H]$ of $X$. Define $C$ as the proper transform of this section. Then $C$ is irreducible. We have $C\in[H'-E]$ and $C^2=w-1<0$. By \cite[Lem. 1.22]{kollar2008}, $\NEb(X')$ has extremal rays $\R_{\geq0}[C]$ and $\R_{\geq0}[E]$. $\Nef(X')$ is dual to $\NEb(X')$ under the pairing given by the intersection product, hence, spanned by the extremal rays $\R_{\geq0}[H'-wE]$ and $\R_{\geq0}[H']$.

By the Theorem of Zariski-Fujita \cite[Rem. 2.1.32]{lazarsfeldpos}, a movable divisor on a normal projective surface is semiample.
Hence, $X'$ is a MDS if and only if the ray $\R_{\geq0}[H'-wE]$ contains a semiample divisor,
which is further equivalent to the following: for some $k>0$, some effective Cartier divisor $F\in [kH'-kwE]$ does not contain $C$ as a component. Note that for any integer $r>0$, such $F$ exists if and only if some $G\in [krH'-krwE]$ exists which does not contain $C$ as a component.
Hence we can replace the triangle $\Delta_0$ by $r\Delta_0$. In particular,
recall $\Delta_1=m\Delta_0$ is the smallest good lattice triangle of the given slopes. We will then replace $\Delta_0$ by $\Delta_1$.

An effective divisor $F$ in the class $[km(H'-wE)]$ is defined by the Laurent polynomial
\[
	f(x,y)=\sum_{(i,j)\in k\Delta_1\cap \Z^2} b_{ij} x^i y^j,
\]
whose partial derivatives up to order $kmw-1$ all vanish at $e=(1,1)$. Notice that the curve $C$ passes through the two torus invariant points, which corresponds to the two vertices $kp$ and $kq$ of $k\Delta_1$ different from $(0,0)$. As pointed out in \cite{GK}, $f$ vanishes on $C$ if and only if $f$ vanishes at $kp$ or $kq$, which is equivalent to that the coefficient $b_{kp}$ or $b_{kq}$ is zero. Combining these results, one has the following result:
\begin{lemma}[\cite{GK}]
	\label{gk1}
	The blow-up $X'$ is a MDS if and only if there exists an integer $k>0$, and some effective divisor $F$ in the class $[km(H'-wE)]$ given by $f(x,y)=\sum_{(i,j)\in k\Delta_1\cap \Z^2} b_{ij} x^i y^j$, such that the coefficients $b_{kp}$ and $b_{kq}$ are nonzero. 
\end{lemma}

Further Gonz\'{a}lez and Karu proved that if there exists a derivative $D$ with order up to $kmw-1$, which vanishes at every monomial $x^i y^j$ indexed by $(i,j)\in k\Delta\cap \Z^2$ when evaluated at $(x,y)=(1,1)$, but does not vanish at $(i,j)=P$ or $Q$, then all $f(x,y)$ whose partial derivatives up to order $kmw-1$ will vanish at $e=(1,1)$ must have coefficients  $b_{kp}$ or $b_{kq}$ zero. In fact, this can be translated into the following statements.
\begin{lemma}[See also \cite{castravet16}]
	Fix an integer $k>0$ and a vertex $p\neq (0,0)$ of $\Delta_1$. Then the following two conditions are equivalent:
	\begin{enumerate}[(1)]
		\item  For every
			\[
			f(x,y)=\sum_{(i,j)\in k\Delta_1\cap \Z^2} b_{ij} x^i y^j,\]
			such that the partial derivatives up to order $kmw-1$ all vanish at $e=(1,1)$,
		$f(x,y)$ has zero coefficient $b_{kp}$;
		\item There exists a derivative $D$ of order up to $kmw-1$, 
			\[
				D=\sum_{0\leq u,v\leq kmw-1} \alpha_{u,v} \partial_x^u \partial_y^v,
			\]
			such that $D(x^i y^j)\mid_{(1,1)}=0$ for all $(i,j)\in k\Delta_1\cap \Z^2$ except $kp$, and $D(x^i y^j)\mid_{(1,1)}\neq 0$ at $kp$.
	\end{enumerate}
	\label{iff}
\end{lemma}

\pf. Consider the $\CC$-vector space $U$ of all $f(x,y)=\sum_{(i,j)\in k\Delta_1\cap \Z^2} b_{ij} x^i y^j$, and the $\CC$-vector space $V$ of all derivatives $	D=\sum_{0\leq u,v\leq kmw-1} \alpha_{u,v} \partial_x^u \partial_y^v$. There is a natural pairing
$\langle\cdot,\cdot\rangle : U\times V\ra \CC$ by $\vr{f,D}= D(f)\mid_{(x,y)=(1,1)}$, which induces natural morphisms $\phi: V\ra U^*$ and $\tau: U\ra V^*$. Let $W\subset U$ be the codimension-one subspace of $U$ spanned by $x^i y^j$ with $(i,j)\in k\Delta_1\cap \Z^2$ except $kp$. Define $W^0:=\{D\in V\mid \vr{f,D}=0, \txt{for all} f\in W\}$.
Then condition (1) in Proposition \ref{iff} is equivalent to $\ker \tau\subseteq W$, and condition (2) is equivalent to that $\ker\phi\subsetneqq W^0$. So we need only show $\ker \tau\subseteq W$ if and only if $\ker\phi\subsetneqq W^0$.
 
Indeed, if there exists a $D_0\in W^0$ but $D_0\not\in \ker \phi$. Since $W$ is codimension-one in $U$, $U=W\oplus \CC x$ for some $x\in U-W$, then $\vr{x,D_0}\neq 0$. Now for every $f\in\ker \tau$, $f=zx+g$ for some $z\in\CC, g\in W$. Then $\vr{zx,D_0}=\vr{f,D_0}-\vr{g,D_0}=0$, hence $z=0$, so $f\in W$. Therefore $\ker \tau\subseteq W$.
Conversely, suppose $\ker \tau\subseteq W$, then $\tau(W)\cong W/\ker \tau$ is a subspace in $V^*$, which is a proper subspace of $\tau(U)\cong U/\ker \tau$. Let $\vr{\cdot,\cdot}_v$ be the dual paring between $V$ and $V^*$. Then $\vr{\cdot,\cdot}_v$ is perfect. Let $\tau(U)^0=\{D\in V\mid \vr{D,d}_v=0, \txt{for all} d\in \tau(U)\}$, and $\tau(W)^0=\{D\in V\mid \vr{D,d}_v=0, \txt{for all} d\in \tau(W)\}$. Then we find $\tau(U)^0\subsetneqq \tau(W)^0=W^0$. Since $\ker \phi\subseteq \tau(U)^0$, we conclude that $\ker\phi\subsetneqq \tau(W)^0$.\qed 
 
\begin{prop}
	\label{gk2}
	The blow-up $X'$ is not a MDS if and only if for every integer $k>0$, there exists a nonzero vertex $p$ of $\Delta_1$ and a derivative $D$ of order up to $kmw-1$:
	\[
	D=\sum_{0\leq u,v\leq kmw-1} \alpha_{u,v} \partial_x^u \partial_y^v,
\]
	such that $D(x^i y^j)\mid_{(1,1)}=0$ for all $(i,j)\in k\Delta_1\cap \Z^2$ except $kp$, and $D(x^i y^j)\mid_{(1,1)}\neq 0$ at $kp$.
\end{prop}
\pf. The sufficiency  is clear from Lemma \ref{gk1} and Proposition \ref{iff}.
For the necessity, suppose $X'$ is not a MDS, then there exists $k>0$, such that every such $f(x,y)$ has zero coefficients $b_{kp}$ or $b_{kq}$. Because all such $f(x,y)$ form a vector space, it must be the case that either  $b_{kp}=0$ for all $f(x,y)$ or $b_{kq}=0$ for all $f(x,y)$. Without loss of generality, we assume $b_{kp}=0$ for all $f(x,y)$. Then Lemma \ref{iff} implies the existence of such derivative $D$, hence the necessity follows. \qed

Finally, let $M$ be the matrix associated to the paring $\vr{\cdot,\cdot}$ in the proof of Proposition \ref{iff}. 
If we adopt the notation that $(a)_b:=a(a-1)\cdots (a-b+1)$, then $M_{(i,j),(u,v)}=(i)_u (j)_v$. (See Section \ref{dayu2})

\begin{lemma}
	\label{column1}
	For every $n> 0$, the vector space spanned by polynomials $\{(x)_u \cdot (y)_v\mid 0\leq u+v\leq n\}$ has a basis $\{x^u y^v\mid 0\leq u+v\leq n\}$.
\end{lemma}
\pf. We need only show every $x^u y^v$ with $u+v\leq n$ is generated by $\{(x)_u (y)_v\mid 0\leq u+v\leq n\}$. We make inductions on $k=u+v$. For $k=0$ the lemma is already true. Suppose this claim holds for $1,\cdots,k$. For $s+t=k+1$, $x^s y^t-(x)_s (y)_t$ is a polynomial of degree at most $k$, and hence is spanned by $\{(x)_u (y)_v\mid 0\leq u+v\leq k\}$ by the induction hypothesis. Hence $x^s y^t$ is spanned by  $\{(x)_u (y)_v\mid 0\leq u+v\leq k+1\}$. \qed

\pf{\em\,of Proposition \ref{intpblm}}. By Lemma \ref{column1}, the existence of such derivative $D$ in Proposition \ref{iff}
is equivalent to the existence of a polynomial 
\[
f=\sum_{u+v\leq kmw-1} \alpha_{u,v} x^u y^v,
\]
which vanishes at every  $(i,j)\in k\Delta_1\cap \Z^2$ except the vertex $kp$. This polynomial $f$ defines a curve $Y$ in $\PP^2 $ with degree up to $kmw-1$, such that $Y$ passes through the complement of $kp$ in $k\Delta_1\cap \Z^2$ but does not pass $kp$.\qed

\section{The Reduced Degree}\label{rd}

We prove Theorem \ref{unity} and \ref{main01} in this section.
For rational slopes $s_1<s_2<s_3$, recall the {\em width}\;of $s_1,s_2,s_3$ equals $w(s_1,s_2,s_3)=1/(s_2-s_1)+1/(s_3-s_2)$, and we assume $w<1$. 
As defined in Section \ref{mainresult}, $\Delta_0$ is the smallest triangle with the given slopes, such that if one vertex is placed at $(0,0)$, then the opposite side passes through $(0,1)$. Next, $\Delta_1$ is the smallest good lattice triangle with the given slopes. Finally, let $m$ be the integer such that $\Delta_1=m\Delta_0$.

Let the leftmost and rightmost vertices of $\Delta_1$ be $p=(x_1,y_1)$ and $q=(x_2,y_2)$ respectively. Then it is easy to show that
\[x_1=\frac{1}{s_1-s_2},\quad x_2=\frac{1}{s_3-s_2},\quad y_1=\frac{s_1}{s_1-s_2},\quad y_2=\frac{s_3}{s_3-s_2}.\]

Recall that the reduced degree of $s_1,s_2,s_3$ equals to the largest nonnegative integer $d$, such that the number of columns in $\Delta_1$ containing at most $d$ points
equals to $d$ (Definition \ref{redd})

We define $l_k=\fl{k s_2}-\cl{k s_1}+1$, and $ r_k=\fl{k s_3}-\cl{k s_2}+1$, for $k=1,2,\cdots$.
\begin{defi}\label{pi}
	$\pi(n)$ is the total number of entries in the two sequences $\{l_k\}_{k\geq 1}$ and $\{r_k\}_{k\geq 1}$ which do not exceed $n$.
\end{defi}

From a geometrical view, $l_k$ and $r_k$ are the number of lattice points on the $k$-th column from the left and right of a lattice triangle $k'\Delta_1$ for all $k'$ such that $k'\Delta_1$ has as least $k$ columns on left and right of the line $\{x=0\}$. It follows that $\pi(n)$  is the number of columns with no more than $n$ points in a sufficiently large lattice triangle $k'\Delta_1$.
Finally, it is clear that $\pi(0)=0$.

\begin{prop}
	Given any rational numbers $s_1<s_2<s_3$ such that the width $w<1$. Assume $\Delta_1=m\Delta_0$. We have
	\begin{enumerate}
		\item $\pi(n)\leq n$ for all $n=0,1,2,\cdots$.
		\item If $n\geq mw$, then  $\pi(n)<n$.
	\end{enumerate}
	\label{period}
\end{prop}

\pf. See Section \ref{2}.

\pf {\em\;of Theorem \ref{unity}.} 
For any $k>0$ fixed, we run the reduction for the triangle $k\Delta_1$ by B\'{e}zout's Theorem as described in Definition \ref{redd}. By Proposition \ref{period}, $\Delta_1=m\Delta_0$, so the width of $k\Delta_1$ equals to $W:=kmw$. Let the reduced degree of $k\Delta_1$ be $d_k$.

By Definition \ref{redd2}, we denote the reduced degree of $s_1,s_2,s_3$ as $d$.
We need only show that all $d_k=d$ for $k=1,2,3,\cdots$.

By sorting the numbers of lattice points in each column of $k\Delta_1$ increasingly, we obtain a sequence $U=\{u_i\}_{i=1}^{W-1}$ of length $W-1=kmw-1$ (because we exclude the two vertices $kp$ and $kq$) and the last term of $U$ equals to $km+1$. Hence at the $i$-th step, we are comparing $u_{W-i}$ with $W-i$. So $d_k$ equals to the largest integer $d_k<km+1$ such that $u_{d_k}\leq d_k$, or equivalently, the largest integer $d_k$ such that there exist at least $d_k$ columns of $\leq d_k$ integer points in $k\Delta_1$. Notice the number of columns in $k\Delta_1$ of at most $s$ integer points is at most $\pi(s)$, which is further bounded by $s$ by Proposition \ref{period}. Hence there can exist at most $d_k$ columns of $\leq d_k$ integer points in $k\Delta_1$. Hence $d_k$ is the largest integer $d_k$ such that there exist exactly $d_k$ columns of $\leq d_k$ integer points in $\Delta_1$, which by Definition \ref{redd2} is the reduced degree $d$ of $s_1,s_2,s_3$. 

Recall the definition of $\pi(n)$ again, we find that $d=\pi(d)$. By Proposition \ref{pinw}, we have $d=\pi(d)<w(d+1)$. Since $w<1$, it follows that $d<w/(1-w)$.
\qed

\begin{corollary}
	Given rational slopes $s_1<s_2<s_3$ with width $w<1$.  Denote the two nonzero vertices of the smallest lattice triangle $\Delta_1$ as $p$ and $q$. Let $d$ be the reduced degree of $s_1,s_2,s_3$. Let $X$ be the toric surface determined by $\Delta_1$.
Then the blow-up $X'=\bl_e X$ at the torus identity point $e$ is not a MDS if and only if for every integer $k>0$ there exists a curve $Y$ in $\PP^2$, of degree up to $d$, and a nonzero vertex $p$ of $\Delta_1$, such that 
\begin{enumerate}
	\item $Y$ passes through all the points in those columns of $k\Delta_1$ of at most $d$ points, and the vertex $kq$, but
	\item $Y$ does not pass through $kp$. 
\end{enumerate}
\label{int}
\end{corollary}

\pf. Suppose the curve $Y$ of degree $\leq W-1=kmw-1$ in Theorem \ref{intpblm} exists. Then we can run the reduction via B\'{e}zout's Theorem and conclude that there exists a curve $Y_d$ of degree $\leq d$ satisfying the conditions in the corollary, when $d$ is the reduced degree. 

Conversely, given a curve $Y_d$ of degree $\leq d$, passing through all points in columns with at most $d$ points, and the vertex $kq$, but not $kp$. Then there are exactly $W-d-1$ columns left in $k\Delta_1$. As a result, the union of $Y_d$ with all the lines lying under every such column is a curve of degree $\leq W-1$, passing through all lattice points in $k\Delta_1$ but not $kp$. Therefore the theorem holds.\qed

\pf{\em\;of Theorem \ref{main01}}.
When $d=0$, such curve $Y$ in Corollary \ref{int} does not exist, hence $X'$ is a MDS.

For (2), suppose $d'=1$. We claim $s_2\not\in \Z$. Indeed, if $s_2\in \Z$, without loss of generality we can assume $s_2=0$. Since $w<1$, we find $s_1<-1$ and $s_3>1$. Therefore $l_1\geq 2$ and $r_1\geq 2$. By Lemma \ref{clm}, $\{l_k\}$ and $\{r_k\}$ are both increasing. So $\pi(1)=0$, which contradicts to $d'=1$. 

Now $d'=1$ implies $l_1=1$ or $r_1=1$. By symmetry we need only prove for the case $l_1=1$. Now In $k\Delta_1$, the first column from the left contains only one point $t$. Since $s_2\not\in \Z$, $t$ is not on the side between $kq$ and $kp$. That is, $kp$, $kq$ and $t$ are not collinear. Then there exists a curve $Y$ whose irreducible components are
\begin{itemize}
	\item The line $L$ through $kq$ and $t$; 
	\item All vertical lines lying under the rest columns of at most $d$ points, except the first column from the left,
\end{itemize}
where $d$ is the reduced degree. 
Then this curve $Y$ has degree $1+(d-1)\cdot 1=d$. Further $Y$ does not passes through $kp$ since $L$ does not passes through $kp$ and  none of the other irreducible components pass through $kp$. Hence $X'$ is not a MDS by Corollary \ref{int}.

Finally if $d=1$, then $d'\leq d=1$, and $d'\neq 0$. So $d'=1$. So it follows from (2). \qed

\section{Examples from Weighted Projective Planes}\label{wpp}
In this section we apply Theorems \ref{main01} and Proposition \ref{main02} to blow-ups of weighted projective planes. 

Throughout this section, we assume that $a,b,c$ are pairwise coprime, such that there exists a relation $(e,f,-g)$ of $(a,b,c)$ (See Application \ref{appwpp}). That is, $ae+bf-cg=0$, with $e,f,g$ positive integers and $\gcd(e,f,g)=1$.

\begin{prop}\label{slopeformula}
	For every pairwise coprime triple $(a,b,c)$ with a relation $(e,f,-g)$, there exists a unique integer $r$ such that $1\leq r\leq g$, $g\mid er-b$ and $g\mid fr+a$. Let $\Delta$ be the triangle with slopes 
	\[
	s_1=\frac{er-b}{eg}, \quad s_2=\frac{r}{g}, \quad s_3=\frac{fr+a}{fg}.\]
	Then $\PP(a,b,c)$ is isomorphic to the toric variety $X _\Delta$.
\end{prop}
\pf. The integer $r$ is the solution of the system of congruence equations:
\begin{equation}	\begin{cases}
	ex\e  b &\m g\\
	fx\e -a &\m g.
\end{cases}\label{cong}
\end{equation}
Let $\alpha=\gcd(e,g)$, $\beta=\gcd(f,g)$, with $e=\alpha e_0$, $f=\beta f_0$. Since $\gcd(e,f,g)=1$, we find $\gcd(\alpha, \beta)=\gcd(\alpha, f)=\gcd(\beta, e)=1$. Hence $\alpha \beta\mid g$. Further $\alpha \mid cg-ae=bf$. Since $\gcd(\alpha,f)=1$, $\alpha\mid b$. Similarly $\beta\mid a$.  Now let $g=\alpha \beta g_0$, $b=\alpha b_0$, and $a=\beta a_0$. It follows that (\ref{cong}) is equivalent to 
\begin{equation}\begin{cases}
		e_0 x\e  b_0 &\m {\beta g_0}\\
		f_0 x\e -a_0 &\m {\alpha g_0}.
	\end{cases}\label{cong2}\end{equation}
By Chinese reminder theorem, the system (\ref{cong2}) has a solution if and only if $e_0\ik b_0\e -f_0\ik a_0$ $\m {\gcd(\beta g_0, \alpha g_0)}$, where the inverse of $g_0$(or $f_0$) is taken in the multiplicative group $(\Z/\beta g_0\Z)^\times$ (respectively, in $(\Z/\alpha g_0\Z)^\times$). Notice that $\gcd(\beta g_0, \alpha g_0)=g_0$. Then $(e_0\ik b_0+f_0\ik a_0)e_0 f_0=f_0 b_0 +e_0 a_0=c g_0$. Hence $(e_0\ik b_0+f_0\ik a_0)e_0 f_0\e 0\m {g_0}$. Since $\gcd(e_0, g_0)=1$, and $\gcd(f_0, g_0)=1$, we find $(e_0\ik b_0+f_0\ik a_0)\e 0\m {g_0}$. Finally, this solution is unique modulo $\lcm (\alpha g_0, \beta g_0)$, which equals to $g$.

By the definition of $s_i$ in the Proposition, we have
\[\vec{n}_1=(\frac{er-b}{g}, -e),\quad \vec{n}_2=(-r,g), \quad \vec{n}_3=(\frac{fr+a}{g},-f)\]
are normal vectors of the sides of $\Delta$. They satisfy the relation that $a \vec{n}_1+c \vec{n}_2+b \vec{n}_3=\vec{0}$. It remains to show that $\vec{n}_1,\vec{n}_2$ and $\vec{n}_3$ span the lattice $N=\Z^2$, and are all primitive vectors.

In order to generate the lattice $\Z^2$, it suffices to show that $\vec{e}_1$ and $\vec{e}_2$ are linear combinations of $\vec{n}_1,\vec{n}_2$ and $\vec{n}_3$ with integer coefficients. Since $a,b,c$ are assumed as pairwise coprime and $ae\e -bf\m c$, we have $a\ik f\e -b\ik e \m c$ (inverses taken in $(\Z/c\Z)^\times$), so the following system of equation of $y$:
\begin{equation}\begin{cases}
		ay\e f &\m c\\
		by\e -e &\m c.
\end{cases}
\end{equation}
has a unique solution (which we still call $y$) mod $c$.
Therefore, there exist integers $x,z$ such that $xc=ay-f$ and $zc=by+e$. Then direct calculation shows $x\vec{n}_1+y\vec{n}_2+z\vec{n}_3=(1,0)=\vec{e}_1$, using that $ae+bf=cg$.

On the other hand, the following system of equation of $y'$:
\begin{equation}\begin{cases}
		ay'\e (fr+a)/g &\m c\\
		by'\e -(er-b)/g &\m c.
\end{cases}
\end{equation}
has a unique solution $y'$ mod $c$. This follows from that $a(er-b)/g+b(fr+a)/g=c\e 0\m c$. There exist integers $x',z'$ such that $x'c=ay'-(fr+a)/g$ and $z'c=by'+(er-b)/g$. It can be calculated then that $x'\vec{n}_1+y'\vec{n}_2+z'\vec{n}_3=(0,1)=\vec{e}_2$.

It remains to show that $\vec{n}_1,\vec{n}_2$ and $\vec{n}_3$ are primitive vectors in $\Z^2$. Indeed, suppose $d>0$ with $d\mid r$ and $d\mid g$, then $d\mid er-b$ and $d\mid fr+a$. So $d\mid b$ and $d\mid a$, therefore $d=1$. This shows $\gcd(r,g)=1$, so $\vec{n}_2$ is primitive. For $\vec{n}_1$, suppose $t>0$ with $t\mid (er-b)/g$ and $t\mid e$. Then $t\mid er-b$, so $t\mid b$. On the other hand, since $f\cdot(er-b)/g-e\cdot(fr+a)/g=-c$ and $(fr+a)/g$ is an integer, we have $t\mid c$. By assumption, $b$ and $c$ are coprime, so $t=1$, and $\vec{n}_2$ is primitive. The result for $\vec{n}_3$ follows from symmetry.
\qed

Now if the slopes are given by Proposition \ref{slopeformula}, then we have
\[w=\frac{1}{s_2-s_1}+\frac{1}{s_3-s_2}=\frac{eg}{b}+\frac{fg}{a}=\frac{cg^2}{ab}.\]
Therefore, all our results and definitions apply for weighted projective surfaces $\PP(a,b,c)$ such that $w=cg^2/ab<1$.

\begin{corollary}\label{g1}
	For every pairwise coprime triple $(a,b,c)$ with a relation $(e,f,-g)$ such that $g=1$ and $c<ab$, the reduced degree of the corresponding slopes is zero, and the blow-up $X'$ is a MDS.
\end{corollary}
\pf. If $g=1$ and $c<ab$, then the width $w=c/ab<1$. By Proposition \ref{slopeformula}, $r=1$, so $s_2=1$. By Corollary \ref{zs2}, the reduced degree equals zero, and $X'$ is a MDS.\qed

\begin{remark}\label{big}
	The MDS claim of Corollary \ref{g1} also follows from Cutkosky's results. It is shown in \cite{cutkosky1991} that if $-K_{X'}$ is a big divisor, then $X'$ is a MDS. We claim if the width $w<1$, then $-K_{X'}$ is big if and only if $cg<a+b+c$. Therefore if $w<1$ and $cg<a+b+c$, then $X'$ is a MDS. In particular, when $g=1$, $-K_{X'}$ is big, and $X'$ is a MDS.

Indeed, let $A=\oO_X (1)$. Let $H$ be defined as in Section \ref{cox}. Then $H=\alpha A$ for some $\alpha\in \Q$. Since $A^2=1/abc$ \cite[Lem. 9]{cutkosky1991} and $H^2=w=cg^2/ab$, we find $r=cg$. Therefore $H=cgA$. The canonical divisor of $X$ is $K_X\equiv\oO_X(-a-b-c)=-(a+b+c)A=-\frac{a+b+c}{cg}H$. Hence $-K_{X'}\equiv\frac{a+b+c}{cg}H'-E$. 

In Section \ref{cox}, we showed that when the width $w<1$, there exists a negative curve $C$ on $X'$ in the class $H'-E$. Therefore $C$ and $E$ span the two extremal rays of $\overline{\NE}(X')$. Now $-K_{X'}$ is big if and only if $-K_{X'}$ is in the interior of $\overline{\NE}(X')$. That is, $cg<a+b+c$.\qed
\end{remark}

\begin{defi}
	Let $a,f,g,r$ be positive integers such that $g\mid fr+a$. Define $\au(a,f,g,r)$ as the set of pairwise coprime triples $(a,b,c)$ such that there exists a relation $(e,f,g)$ with $w<1$. 

We further use $\au(a,f,g,r)_0$, $\au(a,f,g,r)_1$ and $\au(a,f,g,r)_{\geq 2}$ to denote the subsets of $\au(a,f,g,r)$, consisting of triples of minimal degree $0$, $1$, and at least $2$ respectively.
\end{defi}

It is shown in \cite{cutkosky1991} and \cite{HS} that when one of $a,b,c$ is $\leq 4$ or equal to $6$,  $-K_{X'}$ is big, and $X'=\bl_e \PP(a,b,c)$ is MDS.  So the smallest unknown case is when one of $a,b,c$ is $5$.
We apply Proposition \ref{main02} to classify all $\PP(5,b,c)$ with $w<1$, by whether the reduced degree $d$ is zero, one, or at least $2$.

\begin{example}\label{a5}
Assume $a=5$ (equivalently $b=5$ by symmetry).

An important observation emerges that $fg/a<w=fg/a+eg/b<1$, so that $fg<a=5$. Therefore, there are only finite many choices of $f$ and $g$.
The case $g=2, f=2$ contradicts to the assumption that $\gcd(e,f,g)=1$. 
On the other hand, the case $g=1$ gives reduced degree zero by Corollary \ref{g1}.

So the remaining cases are: 1). $f=1, g=2$; 2). $f=1,g=3$; 3). $f=1, g=4$.

\para {\bf Case I. $f=1, g=2$}.
Here $r=1$, so $s_2=1/2$, and $s_3=3$. In addition, $\gamma=4$.

We can check that (4) of Proposition \ref{main02} is satisfied. So Proposition \ref{main02} shows that the reduced degree is zero if and only if $s_1\leq -6/5$.
Notice when $s_2$ and $s_3$ are fixed, the following are equivalent:
\begin{enumerate}
	\item $s_1\leq -6/5$;
	\item the width $w\leq 84/85$;
	\item  $b/c\geq 17/21$.
\end{enumerate}

Further $\au(5,1,2,1)_1=\emptyset$ since $l_1\geq 2$ and $r_1\geq 2$. As a result, a triple in this case has $d'\geq 2$ if and only $84/85<w<1$, or equivalently $4/5< b/c< 17/21$. Hence the above argument shows
\begin{align*}
	\au(5,1,2,1)_0=&\{(5,b,c)\mid b/c\geq 17/21 \txt{and} 5\mid 2c-b\}.\\
	\au(5,1,2,1)_{\geq 2}=&\{(5,b,c)\mid 4/5< b/c< 17/21 \txt{and} 5\mid 2c-b\}.\\
	\au(5,1,2,1)_1=&\emptyset.
\end{align*}
It is easy to prove by slopes that all the triple in $\au(5,1,2,1)_{\geq 2}$ have minimal degree $8$.

Calculation shows the triples with smallest $b$ and $c$ in  $\au(5,1,2,1)_{\geq 2}$ are:
\[(5,37,46), (5,54,67), (5,57,71), (5,71,88), \dots.\]

\para {\bf Case II. $f=1, g=3$}. We have $r=1$, $s_2=1/3$, $s_3=2$, and $\gamma=9$.

(4) of Proposition \ref{main02} is satisfied. 
The reduced degree is zero if and only if $s_1\leq -11/5$, or equivalently
\begin{enumerate}
	\item the width $w\leq 189/190$;
	\item  $b/c\geq 38/21$.
\end{enumerate}
We have 
\begin{align*}
	\au(5,1,3,1)_0=&\{(5,b,c)\mid b/c\geq 38/21 \txt{and} 5\mid 3c-b\}.\\
	\au(5,1,3,1)_{\geq 2}=&\{(5,b,c)\mid 9/5<b/c< 38/21 \txt{and} 5\mid 3c-b\}.\\
\au(5,1,3,1)_1=&\emptyset.
\end{align*}
Similarly, it is easy to show that all the triples in $\au(5,1,3,1)_{\geq 2}$ have minimal degree $12$.
The smallest examples in $\au(5,1,3,1)_{\geq 2}$ are
\[(5, 83, 46), (5, 121, 67), (5, 128, 71), \cdots.\] 

\para {\bf Case III. $f=1, g=4$}. We have $r=3$, $s_2=3/4$, $s_3=2$, and $\gamma=16$.

In this case, Proposition \ref{main02} (1) shows that all triples have reduced degree zero. 
So $\au(5,1,4,3)=\au(5,1,4,3)_0$.

\para
{\bf Conclusion.} When $5,b,c$ are coprime and $w=cg^2/(5b)<1$, $X'$ is a MDS unless in the following two cases.
\begin{enumerate}[(i)]
	\item $4/5< b/c< 17/2$ and  $5\mid 2c-b$;
	\item $9/5<b/c< 38/21$ and $5\mid 3c-b$.
\end{enumerate}
\end{example}

\begin{example} Suppose (4) of Proposition \ref{main02} holds. It follows from Proposition \ref{main02} that every nonempty $\au(a,f,g,r)_\alpha$ has the form
	\[\au(a,f,g,r)_\alpha=\{(a,b,c)\mid b/c\in I \txt{and} a\mid cg-bf\},\]
where $\alpha\in \{1,\geq 2\}$, and $I$ is an interval in $(g^2/a,\infty)$.

Therefore it suffices to determine the range $I$ of $b/c$ for each set $\au(a,f,g,r)_\alpha$.

For all the possible combinations of $(a,f,g,r)$ such that $a\leq 15$, we find by a computer program that (4) of Proposition \ref{main02} holds. Therefore, we can apply Proposition \ref{main02} when $a\leq 15$. We use a computer program to obtain the following tables (Tables \ref{d1table}, \ref{d2table}) of all nonempty $\au(a,f,g,r)_1$ and $\au(a,f,g,r)_{\geq 2}$, for $a\leq 15$.

In other words, any pairwise coprime triple $(a,b,c)$ such that $a\leq 15$, and $w<1$, which appear in neither of the two tables, gives a blow-up which is MDS.

	\begin{table}[h]
	\centering
	\begin{tabular}{|*{2}{*{2}{|l}|}|}\hline
$(a;f,g;r)$ & range of $b/c$ & $(a;f,g;r)$ & range of $b/c$ \\\hline
$\text{(7; 1, 2; 1)}$ & $\frac{4}{7}<\frac{b}{c}<\frac{3}{5}$ & $\text{(13; 1, 3; 2)}$ & $\frac{9}{13}<\frac{b}{c}<\frac{5}{6}$ \\
$\text{(7; 3, 2; 1)}$ & $\text{All}$ & $\text{(13; 4, 3; 2)}$ & $\text{All}$ \\
$\text{(7; 2, 3; 1)}$ & $\text{All}$ & $\text{(13; 1, 4; 3)}$ & $\frac{16}{13}<\frac{b}{c}<\frac{7}{5}$ \\
$\text{(9; 1, 2; 1)}$ & $\frac{4}{9}<\frac{b}{c}<\frac{1}{2}$ & $\text{(13; 3, 4; 1)}$ & $\text{All}$ \\
$\text{(10; 4, 2; 1)}$ & $\text{All}$ & $\text{(13; 2, 5; 1)}$ & $\text{All}$ \\
$\text{(10; 1, 3; 2)}$ & $\frac{9}{10}<\frac{b}{c}<1$ & $\text{(14; 2, 2; 1)}$ & $\frac{2}{7}<\frac{b}{c}<\frac{3}{10}$ \\
$\text{(10; 2, 4; 1)}$ & $\text{All}$ & $\text{(14; 6, 2; 1)}$ & $\text{All}$ \\
$\text{(11; 1, 2; 1)}$ & $\frac{4}{11}<\frac{b}{c}<\frac{3}{7}$ & $\text{(14; 1, 3; 1)}$ & $\frac{9}{14}<\frac{b}{c}<\frac{2}{3}$ \\
$\text{(11; 5, 2; 1)}$ & $\text{All}$ & $\text{(14; 4, 3; 1)}$ & $\text{All}$ \\
$\text{(11; 2, 5; 2)}$ & $\text{All}$ & $\text{(15; 1, 2; 1)}$ & $\frac{4}{15}<\frac{b}{c}<\frac{1}{3}$ \\
$\text{(12; 3, 3; 1)}$ & $\text{All}$ & $\text{(15; 7, 2; 1)}$ & $\text{All}$ \\
$\text{(13; 1, 2; 1)}$ & $\frac{4}{13}<\frac{b}{c}<\frac{3}{8}$ & $\text{(15; 2, 7; 3)}$ & $\text{All}$ \\
$\text{(13; 5, 2; 1)}$ & $\text{All}$ &   &   \\\hline
\end{tabular}	
\caption{All nonempty $\au(a,f,g,r)_1$ such that $a\leq 15$.
Every pairwise coprime triple $(a,b,c)$ with a relation $(e,f,-g)$ such that $w<1$, which appear in this table, has minimal degree $1$, and gives a blow-up which is not a MDS.
}
	\label{d1table}
\end{table}

\begin{table}[h]\centering
\begin{tabular}{|*{3}{*{2}{|l}|}|}\hline
	$(a; f, g ;r)$ & range of $b/c$ &$(a;f,g;r)$ & range of $b/c$ & $(a;f,g;r)$ & range of $b/c$ \\\hline
$\text{(5; 1, 2; 1)}$ & $\frac{4}{5}<\frac{b}{c}<\frac{17}{21}$ & $\text{(11; 1, 4; 1)}$ & $\frac{16}{11}<\frac{b}{c}<\frac{14}{9}$ & $\text{(14; 4, 2; 1)}$ & $\frac{2}{7}<\frac{b}{c}<\frac{5}{17}$ \\
$\text{(5; 1, 3; 1)}$ & $\frac{9}{5}<\frac{b}{c}<\frac{38}{21}$ & $\text{(11; 1, 5; 4)}$ & $\frac{25}{11}<\frac{b}{c}<\frac{104}{45}$ & $\text{(14; 2, 3; 2)}$ & $\frac{9}{14}<\frac{b}{c}<\frac{19}{29}$ \\
$\text{(7; 1, 3; 2)}$ & $\frac{9}{7}<\frac{b}{c}<\frac{38}{29}$ & $\text{(11; 1, 6; 1)}$ & $\frac{36}{11}<\frac{b}{c}<\frac{149}{45}$ & $\text{(14; 2, 4; 1)}$ & $\frac{8}{7}<\frac{b}{c}<\frac{67}{58}$ \\
$\text{(7; 1, 4; 1)}$ & $\frac{16}{7}<\frac{b}{c}<\frac{67}{29}$ & $\text{(11; 1, 7; 3)}$ & $\frac{49}{11}<\frac{b}{c}<\frac{41}{9}$ & $\text{(14; 2, 4; 3)}$ & $\frac{8}{7}<\frac{b}{c}<\frac{11}{9}$ \\
$\text{(7; 1, 5; 3)}$ & $\frac{25}{7}<\frac{b}{c}<\frac{18}{5}$ & $\text{(11; 1, 8; 5)}$ & $\frac{64}{11}<\frac{b}{c}<\frac{327}{56}$ & $\text{(14; 1, 5; 1)}$ & $\frac{25}{14}<\frac{b}{c}<\frac{17}{9}$ \\
$\text{(8; 2, 2; 1)}$ & $\frac{1}{2}<\frac{b}{c}<\frac{5}{9}$ & $\text{(11; 1, 9; 7)}$ & $\frac{81}{11}<\frac{b}{c}<\frac{52}{7}$ & $\text{(14; 2, 5; 3)}$ & $\frac{25}{14}<\frac{b}{c}<\frac{9}{5}$ \\
$\text{(8; 1, 3; 1)}$ & $\frac{9}{8}<\frac{b}{c}<\frac{11}{9}$ & $\text{(12; 1, 5; 3)}$ & $\frac{25}{12}<\frac{b}{c}<\frac{29}{13}$ & $\text{(14; 1, 9; 4)}$ & $\frac{81}{14}<\frac{b}{c}<\frac{53}{9}$ \\
$\text{(8; 1, 5; 2)}$ & $\frac{25}{8}<\frac{b}{c}<\frac{29}{9}$ & $\text{(12; 1, 7; 2)}$ & $\frac{49}{12}<\frac{b}{c}<\frac{55}{13}$ & $\text{(14; 1, 11; 8)}$ & $\frac{121}{14}<\frac{b}{c}<\frac{26}{3}$ \\
$\text{(9; 1, 4; 3)}$ & $\frac{16}{9}<\frac{b}{c}<\frac{67}{37}$ & $\text{(13; 3, 2; 1)}$ & $\frac{4}{1}<\frac{b}{c}<\frac{5}{14}$ & $\text{(15; 3, 2; 1)}$ & $\frac{4}{15}<\frac{b}{c}<\frac{17}{63}$ \\
$\text{(9; 1, 5; 1)}$ & $\frac{25}{9}<\frac{b}{c}<\frac{104}{37}$ & $\text{(13; 2, 3; 1)}$ & $\frac{9}{13}<\frac{b}{c}<\frac{7}{9}$ & $\text{(15; 3, 3; 1)}$ & $\frac{3}{5}<\frac{b}{c}<\frac{38}{63}$ \\
$\text{(9; 1, 7; 5)}$ & $\frac{49}{9}<\frac{b}{c}<\frac{11}{2}$ & $\text{(13; 1, 5; 2)}$ & $\frac{25}{13}<\frac{b}{c}<\frac{19}{9}$ & $\text{(15; 3, 3; 2)}$ & $\frac{3}{5}<\frac{b}{c}<\frac{11}{16}$ \\
$\text{(10; 2, 2; 1)}$ & $\frac{2}{5}<\frac{b}{c}<\frac{17}{42}$ & $\text{(13; 1, 6; 5)}$ & $\frac{36}{13}<\frac{b}{c}<\frac{149}{53}$ & $\text{(15; 1, 4; 1)}$ & $\frac{16}{15}<\frac{b}{c}<\frac{19}{16}$ \\
$\text{(10; 2, 3; 1)}$ & $\frac{9}{10}<\frac{b}{c}<\frac{19}{21}$ & $\text{(13; 1, 7; 1)}$ & $\frac{49}{13}<\frac{b}{c}<\frac{202}{53}$ & $\text{(15; 1, 7; 6)}$ & $\frac{49}{15}<\frac{b}{c}<\frac{202}{61}$ \\
$\text{(10; 1, 7; 4)}$ & $\frac{49}{10}<\frac{b}{c}<5$ & $\text{(13; 1, 8; 3)}$ & $\frac{64}{13}<\frac{b}{c}<\frac{46}{9}$ & $\text{(15; 1, 8; 1)}$ & $\frac{64}{15}<\frac{b}{c}<\frac{263}{61}$ \\
$\text{(11; 3, 2; 1)}$ & $\frac{4}{11}<\frac{b}{c}<\frac{5}{13}$ & $\text{(13; 1, 9; 5)}$ & $\frac{81}{13}<\frac{b}{c}<\frac{32}{5}$ & $\text{(15; 1, 11; 7)}$ & $\frac{121}{15}<\frac{b}{c}<\frac{131}{16}$ \\
$\text{(11; 1, 3; 1)}$ & $\frac{9}{11}<\frac{b}{c}<\frac{47}{56}$ & $\text{(13; 1, 10; 7)}$ & $\frac{100}{13}<\frac{b}{c}<\frac{47}{6}$ & $\text{(15; 1, 13; 11)}$ & $\frac{169}{15}<\frac{b}{c}<\frac{34}{3}$ \\
$\text{(11; 2, 3; 2)}$ & $\frac{9}{11}<\frac{b}{c}<\frac{8}{9}$ & $\text{(13; 1, 11; 9)}$ & $\frac{121}{13}<\frac{b}{c}<\frac{75}{8}$ &   &   \\\hline
\end{tabular}
\caption{All nonempty $\au(a,f,g,r)_{\geq 2}$ such that $a\leq 15$. 
In other word, every pairwise coprime triple $(a,b,c)$ with a relation $(e,f,-g)$ such that $w<1$ and $a\leq 15$, which do not appear in this table and Table \ref{d1table}, has reduced degree zero, and gives a blow-up which is a MDS.}
\label{d2table}
\end{table}
\label{exdd}
\end{example}

\section{Non-examples when the Minimal Degree $d'\geq 2$}\label{dayu2}
\subsection{The interpolation problem of polynomials of degree $d'$}
\label{enb0}

When the minimal degree $d'\geq 2$, the reduced degree $d\geq 2$ too. Corollary \ref{int} motivates us to search for non-examples of MDS by solving the reduced interpolation problems for every $k>0$. Applying the proof of Theorem \ref{main01} and Theorem \ref{intpblm}, we have a sufficient condition for the blow-up $X'$ not to be a MDS:
\begin{corollary}\label{int2}
	Give rational numbers $s_1<s_2<s_3$ with width $w<1$. Let $X$ the toric variety defined by $s_1,s_2,s_3$. Let $X'$ be the blow-up of $X$ at the torus identity point. 

	Assume the minimal degree $d'\geq 2$. If for every integer $k>0$, there exists a nonzero vertex $p$ of $\Delta_1$ and a curve $Y$ in $\PP^2$, of degree up to $d'$, such that 
\begin{enumerate}
	\item $Y$ passes through all the points in columns of $k\Delta_1$ containing $\leq d'$ lattice points, and the vertex $kq$,
	\item $Y$ does not pass $kp$. 
\end{enumerate}
then the blow-up $X'$ is not a MDS.
\end{corollary}
As a result, we try to solve the interpolation problems in Proposition \ref{int2} in search of non-examples of MDS.
However, in general the reduced degree do not equals to the minimal degree. So the reverse direction of (\ref{int2}) can fail.

Recall Proposition \ref{gk2}. Suppose a curve $Y$ in $\CC^2$ is given by a bivariate polynomial $f(x,y)=\sum_{u+v\leq n}a_{u,v} x^u y^v=0$, of total degree $n$. Let $N={n+1 \choose 2}$. Let $I$ be a set of $N$ distinct points in $\CC^2$. Then we can define an $N\times N$ matrix $M$, whose rows are parametrized by $(i,j)\in I$, and columns parametrized by $J=\{(u,v)\mid 0\leq u+v\leq  n\}$, via
\[M_{(i,j),(u,v)}=i^u\cdot j^v.\]
We say in the following that $M=M_{I,J}$ is the matrix parametrized by $I$ and $J$, where $I,J$ are sets of tuples of the same size.
It follows from linear algebra that $Y$ passes through all points in $I$ if and only $M\xi^T=0$, where $\xi=(a_{u,v})_J$. Further if $\det M\neq 0$, then no curve of degree $\leq n$ passing through all the $l$ points.

In our case, for every $k>0$, let $I_k$ be the set of points in the columns of $\leq d'$ points in $k\Delta_1$, together with the two vertices $kp$ and $kq$. Let $J=\{(u,v)\mid 0\leq u+v\leq d'\}$. We obtain a matrix $M'_k$ parametrized by $I_k$ and $J$. Now if $\det M'_k\neq 0$ for all $k\geq 0$, then there is a unique curve $D$ of degree $\leq d$ passing through $I-\{kp\}$. However, there is no curve of degree $\leq d$ passing through all points in $I$. Hence $D$ does not pass $kq$. By Corollary \ref{int2}, $X'$ is not a MDS.
In summary, we have proved

\begin{corollary}
If $\det M'_k\neq 0$ for all $k\geq 0$, then the blow-up $X'$ is not a MDS. 
\label{int3}
\end{corollary}

\begin{remark}\label{shift}
	Corollary \ref{int3} motivates us to calculate $\det M'_k$. Indeed the minimal degree $d'=s+t$, where $s$ is the number of columns of $\leq d'$ points on the left, and $t$ the number on the right. Recall that a shear translation $(x,y)\mapsto (x,y+lx)$ for $l\in\Z$ on $I$ keeps the property that $M\neq 0$. Hence we shift the triangle $k\Delta_1$ so that 
\begin{enumerate}
	\item the right-most vertex $kq$ is at $(t,0)$;
	\item $0\leq s_3<1/t$, so that the first $t$ columns from the right are all in the first quadrant $\{(x,y)\mid x, y\geq0\}$.
\end{enumerate}

Notice it may be impossible to satisfy (2), but in all the following examples, (2) will be satisfied.
\end{remark}

Now we introduce Dumnicki's notation from \cite{Dumnicki2006}.

\begin{figure}[h]
\centering
\begin{minipage}{0.45\textwidth}
\centering
\includegraphics[width=2cm]{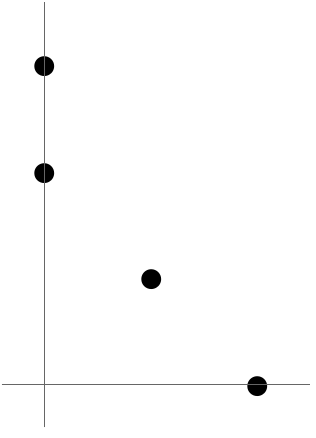} 
\caption{The set $(2^{\uparrow 2}, 1^{\uparrow 1},1^{\uparrow 0})$}
\label{f1}
\end{minipage}\hfill
\begin{minipage}{0.45\textwidth}
\centering
\includegraphics[width=3cm]{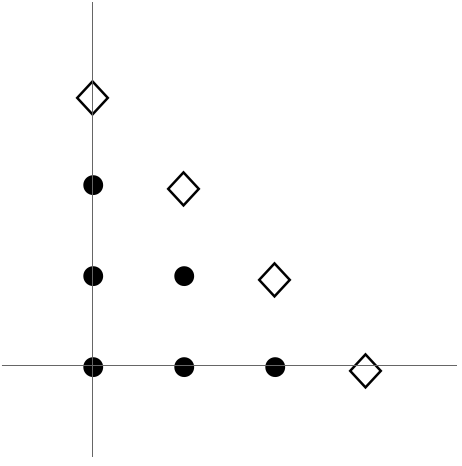} 
\caption{Solid points are $(3^{\uparrow 0}, 2^{\uparrow 0},1^{\uparrow 0}) $. Diamond points index the basis $\ac$}
\label{f2}
\end{minipage}
\end{figure}

\begin{defi}
	Let $a_1,\cdots, a_n, u_1, \cdots, u_n\in\Z$. We define
	\[	(a_1^{\uparrow u_1}, \cdots, a_n^{\uparrow u_n}):=\bigcup_{i=1}^n ( \{i-1\}\times \{u_i, u_i+1,\cdots, u_i+a_i-1\})\subset \Z^2.\]
\end{defi}
For example, the set $(2^{\uparrow 2}, 1^{\uparrow 1},1^{\uparrow 0})$ is shown in Figure \ref{f1}. The set $(3^{\uparrow 0}, 2^{\uparrow 0},1^{\uparrow 0})$ is shown in Figure \ref{f2}.

Here we consider the case when all $u_i=0$. In this case, $(a_1^{\uparrow 0}, \cdots, a_n^{\uparrow 0})$ have two boundaries lying on the two coordinate axis.

Recall that we say a finite set $I$ of planar points {\em imposes independent conditions} on forms of degree $n$ if the linear conditions of $I$ on the coefficients of a polynomial of degree $n$ are independent. Equivalently, for an interpolation problem of a polynomials of degree $n$ vanishing on $I$, each point in $I$ makes the dimension of the solution space drop by exactly one.

Fixing an integer $d>0$, we choose $\sigma=\{a_1,a_2,\cdots, a_n\}$ a subset of $\{1,2,\cdots ,d\}$ containing $d$, and suppose $a_1=d >a_2>\cdots a_n$. We claim that
\begin{lemma}
	Consider the vector space $U$ of bivariate polynomials $f(x,y)$ of degree up to $d$, vanishing on the lattice set $(a_1^{\uparrow 0}, \cdots, a_n^{\uparrow 0})$ defined by $\sigma$.

	Then $U$ have the following basis:
		\[\ac=\Bigl\{ {x\choose u}{y\choose v} \mid (u,v)\in ((d+1)^{\uparrow 0}, \cdots, 1^{\uparrow 0})-(a_1^{\uparrow 0}, \cdots, a_n^{\uparrow 0}))\Bigr\}.\]
		In particular, the lattice set $(a_1^{\uparrow 0}, \cdots, a_n^{\uparrow 0}))$  imposes independent conditions.\label{dual}
\end{lemma}
\begin{remark}
	 If we let $(x)_u=x(x-1)\cdots (x-u+1)$, the basis can also be chosen as
	 \[\ac'=\Bigl\{ (x)_u (y)_v \mid (u,v)\in ((d+1)^{\uparrow 0}, \cdots, 1^{\uparrow 0})-(a_1^{\uparrow 0}, \cdots, a_n^{\uparrow 0}))\Bigr\}.\]
\end{remark}

\pf. Clearly every function in $\ac$ vanishes on every point of $(a_1^{\uparrow 0}, \cdots, a_n^{\uparrow 0})$. It is also clear that they are linearly independent. So it remains to prove that they span the vector space $U$.  The lattice set $S=((d+1)^{\uparrow 0}, \cdots, 1^{\uparrow 0}))$ consists of exactly ${d+1 \choose 2}$ points. By Lemma 4.20 in \cite{Dumnicki2006},
$Y$ imposes independent conditions, so the corresponding interpolation matrix is nonsingular. In particular the subset $(a_1^{\uparrow 0}, \cdots, a_n^{\uparrow 0})$ imposes independent conditions, hence the dimension of $U$ equals to $\vv{\ac}$. 
This shows that $\ac$ span the vector space $U$. \qed

Here we divide the set $I_k$ into union of $I^-_k$ and $I^+_k$, consisting of the points in the left and right corners of $k\Delta_1$ respectively. Now $I^+_k$ satisfies the conditions of Lemma \ref{dual}. Hence $I^+_k$ determines a basis $\ac_k$. 
Now we can define a matrix $M_k$ which is parametrized by $I^-_k$ and $\ac_k$.

\begin{corollary}
	Suppose we can shift the triangles $k\Delta_1$ so that (1) and (2) are satisfied. Let $M_k$ be defined as above. 
If $\det M_k\neq 0$ for all $k\geq 0$, then the blow-up $X'$ is not a MDS.
\label{int4}
\end{corollary}

\pf. Indeed $\ac_k$ span the vector space of bivariate polynomials of degree $\leq d'$ vanishing on $I^+_k$. Hence any polynomial of degree $\leq d'$ vanishing on all $I_k$ is spanned by $\ac_k$. As a result, $\det M'_k=0$ if and only if $\det M_k=0$.\qed

\subsection{Gonz\'{a}lez and Karu's example}
We define the generalized falling factorial $(x)_n$ by letting $(x)_n:={x \choose n}/n!$, for $n\in\Z$. For $n>0$ this coincides with the usual definition $(x)_n=x(x-1)\cdots (x-n)$. For $n=0$, $(x)_0=1$. For $n<0$, we have $(x)_n=1/\left( (x+1)\cdots (x+(-n)) \right)$. So when $n<0$, $(x)_n$ is defined for any $x\not\in \{-1,\cdots,-n\}$.

\label{enb1} In \cite{GK} the non-examples of MDS are given by triangles with minimal degrees $d'=n$, $s=1$ and $t=n-1$. The lattice set $I_k$ is given by
	\[I^-=(1^{\uparrow 0}, n^{\uparrow -n-1})+(A,B),\quad 
	I^+=(n^{\uparrow 0}, \cdots, 1^{\uparrow 0}).\]
where  $A=n-1-kmw$, $B=-s_2 kmw$, and $I_k=I^-\cup I^+$.
Further, $I^+$ defines the basis
\[\ac=\{(x)_n,\cdots,(x)_{n-i} (y)_i,\cdots,(y)_n \}^n_{i=0}.\]
As a result, to prove these examples gives non-MDS blow-ups, we need only show $\det M_k\neq 0$. Instead of working for an individual $k$, we prefer leaving $A$ and $B$ as indeterminate, so that the matrix $M$ parametrized by $I^-$ and $\ac$ is a matrix polynomial in variables $A$ and $B$. Then the determinant of $M$ is a polynomial of $A$ and $B$.

Now the key idea is to partition $I^-$ into $I^-=K_0\cup K_1$. We let $K_0=\{(A,B)\}$, so $K_1=I^-\backslash K_0=\{(A+1,B-2-j)\}_{j=0}^{n-1}$. Let $M_\alpha$ be the submatrix of $M$ parametrized by $K_\alpha$, for $\alpha=0,1$.
\begin{prop}\label{xi}
	\begin{enumerate}
		\item Define a $(n+1)\times 1$ vector
			\[\xi=\bigg((-1)^i {n\choose i} (A+1-n+i)_i (B-2-i)_{n-i}\bigg)_{i=0}^n.\]
			Then $M_1\xi=0$. Furthermore if $A\neq -1,0,1,\cdots, n-2$, $\xi$ is the unique solution of $M_1 \xi=0$ up to a scalar.
		\item 
			$\det M=0$ if and only if $(n+1)(A-n+1)+nB=0$ or $A=-1,0,\cdots,n-2$.
	\end{enumerate}
\end{prop}
\begin{remark}
	Let $A=n-1-kmw$, $B=-s_2 kmw$. Then $(n+1)(A-n+1)-nB)=-kmw(ns_2+n+1)$. Since there are at least $n$ columns in $\Delta_1$, $mw\geq n+1$. Hence $A\leq -2$. Therefore $\det M_k=0$ if and only if $n s_2\in\Z$, which reproves the main theorem in \cite{GK}.
\end{remark}

\pf. The product of the $j$-th row of $M_1$ and $\xi$ is
\begin{align*}
	&\sum_{i=0}^n (A+1)_{n-i} (B-2-j)_{i} (-1)^i {n\choose i} (A+1-n+i)_i (B-2-i)_{n-i}\\
	=&\frac{(A+1)!}{(A+1-n)!}\frac{(B-2-j)!}{(B-2-n)!}\btr{i}\frac{(B-2-i)!}{(B-2-j-i)!}.
\end{align*}
Notice that $(B-2-i)!/(B-2-j-i)!=(B-2-i)_j$ is a polynomial of $i$ with degree $j\leq n-1$. Hence the sum above is zero by Lemma \ref{polyvanish}, which proves that $M_1\xi=0$.

When $A\neq -1,0,\cdots, n-2$, we can divide the $i$-th column of $M_1$ by $(A+1)_{n-i}$ to obtain a matrix $M_2$ of a single indeterminate $B$. That is, $(M_2)_{j,i}=
(B-2-j)_i$. By Lemma \ref{van-stir}, $M_2$ has rank $n$, so $M_1$ has rank $n$, too. Therefore $\xi$ span the space of solutions of $M_1 x=0$.

For (2), when $A\in\{-1,0,1,\cdots,n-2\}$, $(A)_n=(A+1)_n=0$, so the first column of $A$ is zero. When $A=-1$, $(A+1)_i=0$ for all $0\leq i\leq n$, hence $M_1=0$. In both cases, $\det M=0$. When $A\not\in\{-1,0,\cdots,n-2\}$, by (1), $\det M=0$ if and only if $M\xi=0$, which is equivalent to $M_0\xi=0$. Here 
\begin{align*}
	M_0\xi=&\sum_{i=0}^n (A)_{n-i} (B)_{i} (-1)^i {n\choose i} (A+1-n+i)_i (B-2-i)_{n-i}\\
	=&\frac{A!}{(A+1-n)!}\frac{B!}{(B-2-n)!}\btr{i}\frac{(A+1-n+i)}{(B-i-1)(B-i)}\\
	=&(A)_{n-1}(B)_{n+2}\btr{i}\bigg(-\frac{B-A+n-2}{B-1-i}+\frac{B-A+n-1}{B-i}\bigg).\\
\end{align*}
By Lemma \ref{logsum}, 
\[\btr{i}\frac{1}{B-1-i}=-n!(-B)_{-n-1}, \quad \btr{i}\frac{1}{B-i}=-n!(-B-1)_{-n-1}, \]
Hence
\begin{align*}
	M_0\xi=&n!(A)_{n-1}(B)_{n+2}\bigg((B-A+n-2)(-B)_{-n-1}-(B-A+n-1)(-B-1)_{-n-1}\bigg)\\
	=&(-1)^{n+1}n! (A)_{n-1}\bigg((B-A+n-2)(B)-(B-A+n-1)(B-n-1))\bigg)\\
	=&(-1)^{n+1}n! (A)_{n-1}((n+1)(A-n+1)+nB).
\end{align*} \qed

\begin{lemma} Given $n$ distinct number  $a_1,\cdots,a_n\in\R$,
	the matrix $U=(u_{i,j})$ where $u_{i,j}=(a_i)_{j-1}$, $1\leq i,j\leq n$ is nonsingular.
	\label{van-stir}
\end{lemma}
\pf. Recall the identity
\[\sum_{i=0}^n  \stirling{n}{i}(x)_i=x^n\]
for $n>0$, $(x)_i$ the $i$-th falling factorial, and $\stirling{n}{i}$ the Stirling numbers of the second kind, which is the number of partitions of a set of $n$ elements into $i$ nonempty subsets. Hence the Stirling numbers $\stirling{j}{i}$, $0\leq i,j\leq n-1$ give a matrix $P$ which transforms $U$ into a Vandermonde Matrix $V=(v_{i,j})$ where $v_{i,j}=(a_i)^{j-1}$. That is, $UP=V$. The matrix $C$ is nonsingular because all $a_i$s are distinct. The matrix $P$ is unit triangular because $\stirling{j}{i}=0$ when $i>j$ and $\stirling{i}{i}=1$, therefore nonsingular. Hence $U$ is nonsingular. \qed 

\begin{lemma}(See \cite[Lem. 2.3]{GK})
	Let $n>0$ be an integer and $p(x)$ be a polynomial of degree $<n$. Then
	\[\btr{i} p(i)=0.\]
	\label{polyvanish}
\end{lemma}

\begin{lemma}
	\label{logsum}
	Let $n\in\Z_{\geq 0}$. Then
	\[\btr{i} \frac{1}{x+i}=n!(x-1)_{-n-1}=\frac{n!}{x(x+1)\cdots(x+n)}.\]
\end{lemma}
\pf. We make induction on $n$. When $n=0$ this obviously holds.
Suppose the identity holds for $n-1$. Then for $n$, we use \[{n \choose i}={n-1 \choose i}+{n-1 \choose i-1}.\] So
\begin{align*}
	\btr{i}\frac{1}{x+i}&=\sum_{i=0}^{n} (-1)^i {n-1\choose i}\frac{1}{x+i}+\sum_{i=0}^n (-1)^i {n-1\choose i-1}\frac{1}{x+i}\\
	&=\sum_{i=0}^{n-1} (-1)^i {n-1\choose i}\frac{1}{x+i}+\sum_{i=0}^{n-1} (-1)^{i+1} {n-1\choose i}\frac{1}{x+1+i}\\
	&=(n-1)!(x-1)_{-n}-(n-1)!(x)_{-n}\\
	&=n!(x-1)_{-n-1}.
\end{align*}\qed

\subsection{The Case $d'=5$ (Proof of Theorem \ref{main04}, First case)}
\label{enb2}
Here we prove that the first case in Theorem \ref{main04} of $d'=5$ gives non-MDS blow-ups if $5s_2\not\in\Z$.
For the rest two cases $d'=7$ and $9$, we managed to calculate $\det M$ via a computer program. The whole proof of Theorem \ref{main03} and \ref{main04} is given in Section \ref{03pf}.

\begin{figure}[h]
\centering
\begin{minipage}{1\textwidth}
\centering
\includegraphics[width=4cm]{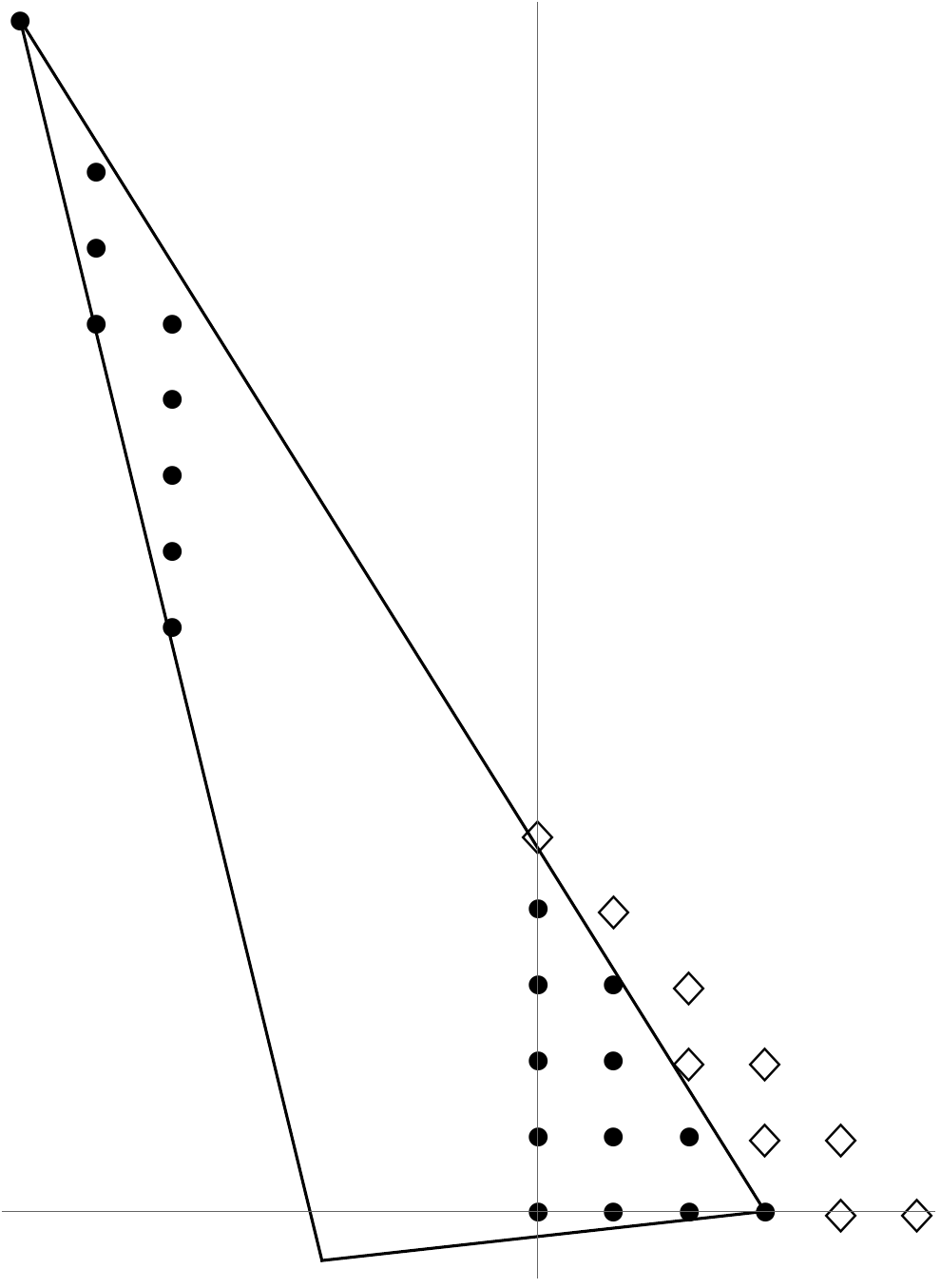} 
\caption{The triangle obtained by adding $-2$ to all the slopes of $k\Delta_1$ and then applying a shear translation. Solid points are the set $I_k$. Diamond points index the basis $\ac$}
\label{f3}
\end{minipage}
\end{figure}

We add $-2$ to all the slopes and translate the triangle so that the assumptions (1) and (2) in Remark \ref{shift} are satisfied.
Define
\[	I^-=(1^{\uparrow 0}, 3^{\uparrow -4},5^{\uparrow -8})+(A,B), I^+=(5^{\uparrow 0},4^{\uparrow 0},2^{\uparrow 0}, 1^{\uparrow 0}).\]
where $-A=kmw-3$ and $B=-s_2 kmw$.
Then in the triangle $k\Delta_1$, the lattice sets $I_k=I^-\cup I^+$, for every $k$. Further $I^+$ define the basis
\[\ac=\{(x)_5,\cdots,(x)_{5-i} (y)_i,\cdots,(y)_5, (x)_4, (x)_3 y, (x)_2 (y)_2\}^5_{i=0}.\]
See Figure \ref{f3}.
In the following we treat $A$ and $B$ as indeterminate. 

The matrix $M$ is defined as the matrix parametrized by $I^-$ and $\ac$. For convenience, we choose the order on $I_k$ as follows:
the vertex $(A,B)$ is the first. Then the first column next to $(A,B)$ from top to bottom. Finally the second column from top to bottom.

We divide the matrix $M$ into blocks:
\[M= \left(
    \begin{array}{c|c}
      E & P\\\hline
      F & Q\\\hline
      G & R
    \end{array}
    \right)
\]
so that $E$ is of size $1\times 6$, $F$ of $3\times 6$, and $G$ of $5\times 6$. That is, $E, F, G$  are indexed by the points with $x$-coordinates $A,A+1,A+2$, respectively. We claim:
\begin{lemma}$ $ 
Define a $6\times 1$ vector
	\[\xi=\bigg((-1)^i {5\choose i} (A+2-5+i)_i (B-4-i)_{n-i}\bigg)_{i=0}^5\]
Then $\det M=0$ if and only if $A=-2,-1,0,1,2$ or $(E-2PQ\ik F)\xi=0$.
	\label{d5lemmagroup}
\end{lemma}
\pf. The same proof of Proposition \ref{xi} shows that $\xi$ is the unique solution of $G\xi=0$ up to a scalar if $(A+2)_5\neq 0$. Solve $(A+2)_5=0$ and we have $A=-2,-1,0,1,2$.

When $A=0,1,2$, the first column of $M$ is zero. When $A=-1$, $M$ has a $8\times 7$ submatrix which is zero. When $A=-2$, $M$ has a $5\times 8 $ submatrix which is zero. Hence in all these cases $\det M=0$. 

So we assume $(A+2)_5\neq 0$. By Lemma \ref{van-stir}, $Q$ is nonsingular. Let $Q\ik F\xi=(a_0,a_1,a_2)^T$. Define the vector $\eta$ such that 
\[\eta^T:=(\xi^T,0,0,0)+(a_0,a_1,a_2,0,0,0,-(A-2)a_0,-(A-1)a_1,-A a_2).\]
We claim that $(F\;Q)\eta=(G \; R)\eta=0$. Let $A_{(i)}$ denote the $i$-th row of a matrix $A$.
Indeed, $(G\; R)_{(i)}\eta=(G\xi)_{(i)}+a_0 (A+2)_5+a_1 (A+2)_4 (B-4-i) +a_2 (A+2)_3 (B-4-i)_2-a_0(A-2) (A+2)_4-a_1(A-1) (A+2)_3 (B-4-i) -a_2(A) (A+2)_2 (B-4-i)_2=0$.

Similarly $(F\; Q)_{(i)}\eta=(F\xi)_{(i)}+a_0 (A+1)_5+a_1 (A+1)_4 (B-2-i) +a_2 (A+1)_3 (B-2-i)_2-a_0(A-2) (A+1)_4-a_1(A-1) (A+1)_3 (B-2-i) -a_2(A) (A+1)_2 (B-2-i)_2$, which equals to $(F\xi)_{(i)}-a_0 (A+1)_4-a_1(A+1)_3 (B-2-i) -a_2(A+1)_2 (B-2-i)_2$. In matrix form, this says $(F\;Q)\eta=F\xi-Q\eta=F\xi-Q(Q\ik F\xi)=0$.

We now claim that the rank of the submatrix $N:=
\begin{pmatrix}
	F & Q\\
	G & R
\end{pmatrix}$
equals to $8$, so $\eta$ spans the solution space of $N$. As a result,
$\det M=0$ if and only if $(E \;P)\eta=0$. Now similar calculation shows that
\begin{align*}
	\begin{pmatrix}
		E & P
	\end{pmatrix}
	\eta=&E\xi-2a_0 (A)_4-2a_1(A)_3 (B) -2a_2(A)_2 (B)_2\\
	=&E\xi-2P(Q\ik F\xi).
\end{align*}

So we prove that $\rk N=8$, assuming $(A+2)_5\neq 0$. We divide the first three columns of $N$ by $-(A-2),-(A-1),-A$ respectively, and add to the $7,8,9$-th column. This reduces $R$ to zero. To find what $Q$ is transformed to, we have
\[	\frac{(A+1)_{5-k} (B-s)_{k}}{-(A-2+k)}+(A+1)_{4-k} (B-s)_{k}
=(A+1)_{3-k}(B-s)_{k},\]
for $s=2,3,4$. 
Therefore $N\sim N_1=\begin{pmatrix}
	F & Q_1\\
	G & 0
\end{pmatrix}$
where \[	Q_1=\begin{pmatrix}\displaystyle 
		(A+1)_3 (B-2)_0 & (A+1)_2 (B-2)_1&(A+1)_1 (B-2)_2\\  
		(A+1)_3 (B-3)_0 & (A+1)_2 (B-3)_1&(A+1)_1 (B-3)_2\\  
		(A+1)_3 (B-4)_0 & (A+1)_2 (B-4)_1&(A+1)_1 (B-4)_2\\  
\end{pmatrix}.\]
Since $(A+2)_5\neq 0$, we can divide on each column of $N$ by  the corresponding factors of $A$ appearing in $G$ and $R_1$. This reduced $G$ and $Q_1$ to Vandermonde-like matrices as in Lemma \ref{van-stir}. Therefore, Lemma \ref{van-stir} shows that $Q$ and the first five columns of $G_1$ are nonsingular. Hence we can reduce the matrix $N_1$ further to the row echelon form, which shows that $N_1$ has full row rank, which equals $8$. Therefore $\rk N=\rk N_1=8$.
\qed

\pf{\em\;of the first case of Theorem \ref{main04}}. By Corollary \ref{int4}, we need to show $\det M\neq 0$ when $5s_2\not\in\Z$. By Lemma \ref{d5lemmagroup}, it suffices to calculate $E\xi-2P(Q\ik F\xi)$. We will find $PQ\ik$, $F\xi$ and $E\xi$ respectively.

Suppose $PQ\ik=\zeta$, then $P^T=Q^T\zeta^T$, we divide the $j$-th column of $Q$ by $(A+1)_{4-j}$ to obtain a matrix $Q_1$ where the $(i,j)$-th entry is $(B-2-i)_j$, for $0\leq i,j\leq 2$. We divide the $i$-th entry of $P$ by $(A+1)_{4-j}$ to get $P_1$. Now $Q_1^T\zeta^T=P_1^T$, and $Q_1^T$ satisfies the assumption of Lemma \ref{van-stir}. Observe that the $j$-th entry of $P_1$ is
\[\frac{(A)_{4-j}}{(A+1)_{4-j}}(B)_j=\frac{A-3+j}{A+1}(B)_j=\frac{1}{A+1}( (A-3)(B)_j+B(B)_j-B(B-1)_j).\]

Hence by linearity and Lemma \ref{d5lemmagroup}, we find
\begin{align*}
(PQ\ik)^T=\zeta^T=&\bigg(\frac{(-1)^j}{2!(A+1)} {2\choose j}\big( \frac{(A-3+B)(4)_3}{2+j}-\frac{B(3)_3}{1+j}\big)\bigg)_{j=0}^2\\
=&\frac{1}{A+1}\bigg(6A+3B-18,-8A-5B+24,3A+2B-9 \bigg).
\end{align*}

Next the $j$-th entry of $F\xi$ equals to
\begin{align*}
	 &\sum_{i=0}^5 (A+1)_{5-i} (B-2-j)_{i} (-1)^i {5\choose i} (A-3+i)_i (B-4-i)_{5-i}\\
	 =&\frac{(A+1)!}{(A-3)!}\frac{(B-2-j)!}{(B-9)!}\wtr{i}\frac{(A-3+i)(B-4-i)!}{(B-2-j-i)!}\\
	 =&(A+1)_{4}(B-j)_{7-j}\wtr{i}(A-3+i)(B-4-i)_{j-2}.
\end{align*}
For $j=0,1,2$.
Using Lemma \ref{stirling}, the above equals to
\begin{align*}
&(A+1)_{4}(B-j)_{7-j}(-1)(B-9)_{j-7}\bigg( (A-3)(6-j)_5+(5-j)_5 5(B-2-j)\bigg).\\
=&-(A+1)_4 (5-j)_4 \bigg( (6-j)(A-3)+5(B-2-j)\bigg).
\end{align*}
In other words, $F\xi=-5!(A+1)_4\bigg( 6A+5B-28, A+B-6,0\bigg).$
For $E\xi$:
\begin{align*}
	 E\xi=&\sum_{i=0}^5 (A)_{5-i} (B)_{i} (-1)^i {5\choose i} (A-3+i)_i (B-4-i)_{5-i}\\
	 =&(A)_3 (B)_{9}\wtr{i}(A-3+i)(A-4-i)(B-4-i)_{-4}.
 \end{align*}
 Use Lemma \ref{stirling} again,
 \begin{align*}
	E\xi =&(A)_3 (B)_{9}\bigg( (A-3)(A-4) (-1)^5(B-9)_{-9}(8)_5 +\\
	 &(2A-7) (-1)^5(B-9)_{-9}(7)_4 (5B)+(-1)^5(B-9)_{-9}(6)_3 (5B)(4B+3)\bigg)\\
	 =&-(A)_3\bigg( (A-3)(A-4) (8)_5 +(7)_4 (5B)+(6)_3 (5B)(4B+3)\bigg).
\end{align*}
Combining the results above, we find
\[E\xi-2 (PQ\ik)(F\xi)=-5!(A)_3 (8(A-3)+5B).\]
Hence $\det M\neq 0$ if $8(A-3)+5B\neq 0$. Let $-A=kmw-3$ and $B=-s_2 kmw$, then this is equivalent to $5s_2+8\neq 0$. Since $-2+1/3<s_2<-2+1/2$, this is equivalent to $5s_2\not\in\Z$. \qed

\begin{remark}
	The determinant $\det M$ can be calculated by a computer program. We used {\em Mathematica 10} \cite{mma104} to obtain that
	 \[\det M=-2^{10}\cdot 3^3\cdot 5 (-2 + A) (-1 + A)^3 A^5 (1 + A)^6 (2 + A)^4 (-24 + 8 A + 5 B),\]
	 which supports our calculation.
\end{remark}

\begin{lemma}
	Let $x\in\R$, $m,n\in\Z_{> 0}$. Define a $(n+1)\times (n+1)$ matrix $U$ by
	\[U_{i,j}:=(x-m-j)_i\]
	where $0\leq i,j\leq n$. Let 
	\[\sigma=\bigg( (-1)^j {n\choose j}\frac{(n+m)_{n+1}}{n!(m+j)}\bigg)_{j=0}^{n}.\]
	Then $U\sigma=(1,(x)_1,\cdots,(x)_{n})^T$.
	\label{van-inv}
\end{lemma}
\pf. The product of the $i$-th row of $U$ with $\sigma$ is
\[\btr{j}(x-m-j)_i\frac{(n+m)_{n+1}}{n!(m+j)}=\frac{(n+m)_{n+1}}{n!}\btr{j}\frac{(x-m-j)_i}{m+j}.\]
Consider as polynomials of $j$, then $(x-m-j)_i=f(j)(m+j)+g(j)$. The quotient $f(j)$ is of degree $i-1$, and $g(j)=(x-m-(-m))_j=(x)_j$ by the reminder theorem. Hence 
\begin{align*}
	\btr{j}\frac{(x-m-j)_i}{m+j}=&\btr{j}f(j)+\btr{j}\frac{(x)_j}{m+i}\\
	=&0+(x)_i\cdot n!(m-1)_{-n-1}\\
	=&\frac{n!(x)_j}{(n+m)_{n+1}},
\end{align*}
where we applied Lemma \ref{polyvanish} and \ref{logsum}. Now the product equals to $(x)_i$.\qed

For the following lemmas, we let $\Delta f(x):= f(x+1)-f(x)$ be the finite difference of $f(x)$ with respect of $x$.
\begin{lemma}
	For $n\neq 0$, $\Delta (x)_n=n (x)_{n-1}$.
	\label{diff}
\end{lemma}
\pf. $\Delta (x)_n=(x+1)_n-(x)_n=((x+1)-(x-n+1))(x)_{n-1}=n(x)_{n-1}$. \qed

\begin{lemma}
	Consider the sum 
	\[g(s):=\btr{i}i^s (x-i)_{-m}\]
	for $n,m>0$, $x\in\R$, and $s\geq 0$.
	\label{stirling}
	Then 
	\begin{enumerate}
		\item $g(0)=(-1)^n (x-n)_{-n-m}(n+m-1)_{n}$;
		\item $g(1)=(-1)^n (x-n)_{-n-m}(n+m-2)_{n-1}n(x+m)$;
		\item $g(2)=(-1)^n (x-n)_{-n-m}(n+m-3)_{n-2}n(x+m)((n-1)x+mn-1)$.
	\end{enumerate}
\end{lemma}
\pf. First for $g(0)$ this is a generalization of Lemma \ref{logsum}.
Indeed, apply Lemma \ref{diff} repeatedly. We have 
\[(x)_{-m}=\frac{\Delta (x)_{-m+1}}{-m+1}=\cdots=\frac{\Delta^{m-1} (x)_{-1} }{(-1)^{m-1}(m-1)!}.\]
Notice that replacing $x$ by $x-i$ for an integer $i$ does not change the finite difference in Lemma \ref{logsum}. Hence we have
\begin{align*}
	g(0)&=\frac{1}{(-1)^{m-1}(m-1)!}\btr{i}\Delta^{m-1} \frac{1}{x-i+1}\\
	&=\frac{1}{(-1)^{m-1}(m-1)!}\Delta^{m-1}\btr{i} \frac{1}{x-i+1}\\
	&=\frac{1}{(-1)^{m-1}(m-1)!}\Delta^{m-1}(n!(-1) (-x-2)_{-n-1})\\
	&=\frac{n!}{(-1)^{m}(m-1)!}\Delta^{m-1}( (-1)^{n+1}(x-n)_{-n-1})\\
	&=\frac{n!}{(-1)^{m}(m-1)!}( (-1)^{n+1}(x-n)_{-n-m})(-1)^{m-1}(n+m-1)_{m-1})\\
	&=(-1)^n (n+m-1)_n (x-n)_{-n-m}.
\end{align*}
For $s=1$, notice $(x-i+1)_{-m+1}=(x-i+1)(x-i)_{-m}$. So $i(x-i)_{-m}=(x+1)(x-i)_{-m}-(x-i+1)_{-m+1}$ for $m\geq 1$. Hence
\begin{align*}
	g(1)=& (x+1)\btr{i}(x-i)_{-m}-\btr{i}(x-i+1)_{-m+1}\\
	=&(x+1)(-1)^n (n+m-1)_n (x-n)_{-n-m}-(-1)^n (n+m-2)_n (x+1-n)_{-n-m+1}\\
	=&(-1)^n (x-n)_{-n-m}(n+m-2)_{n-1}n(x+m).
\end{align*}
Similarly for $s=2$, we have 
\begin{align*}
	i^2(x-i)_{-m}=&i(x+1)(x-i)_{-m}-i(x-i+1)_{-m+1}\\
	=&(x+1)^2(x-i)_{-m}-(2x-3)(x+1-i)_{-m+1}+(x+2-i)_{-m+2}.
\end{align*}
Hence
\begin{align*}
	g(2)=&(x+1)^2(-1)^n (n+m-1)_n (x-n)_{-n-m}\\
	&-(2x-3)(-1)^n (n+m-2)_n (x+1-n)_{-n-m+1}\\
	&+(-1)^n (n+m-3)_n (x+2-n)_{-n-m+2}
\end{align*}
which simplifies to the required result. \qed

\subsection{The Cases $d'=7$ or $9$} We calculated the determinant of the interpolation matrix $M$ by computer programs in the cases when $d'=7$ and $9$, and the slopes are given by the inequalities in Theorem \ref{main04}. The codes for the three cases ($d'=5,7,9$) are available online at \url{https://hezhuangblog.wordpress.com/research}.

In both cases, the construction of the shifted triangle $k\Delta_1$, the index sets $I^-$, $I^+$, $\ac$, and the matrix $M$ are all similar with the case $d'=5$. We used a program in {\em Mathematica 10} to calculate $\det M$. As in the case $d'=5$, we let the coordinate of the vertex in $I^-$ be $(A,B)$. Then $\det M=P(A,B)$ is a polynomial in $A$ and $B$. Let $b_1,\cdots,b_r$ be distinct integers for $r$ sufficiently large. For each $b_i$, evaluating $P(A,b_i)$ at sufficiently many values of $A$ gives a polynomial $Q_{b_i} (A)$ such that $P(A,b_i)=Q_{b_i} (A)$ is the Lagrange polynomial in $A$. Then $P(A,B)$ equals the Lagrange polynomial of the data set $\{b_i,Q_{b_i} (A)\}$. Specifically:

(1) for the case $d'=7$ (within the given constraints on slopes),
\begin{align*}
	\det M=&-2^{26}\cdot 3^{10}\cdot 5^3\cdot 7(-3 + A) (-2 + A)^3 (-1 + A)^6\cdot\\
	&A^9 (1 + A)^{11} (2 + A)^{10} (3 + A)^6 (12 A-48+7B).
\end{align*}

(2) For the case $d'=9$ (within the given constraints on slopes),
\begin{align*}
	\det M=& -2^{59}\cdot 3^{24}\cdot 5^7\cdot 7^3 ((-4 + A) (-3 + A)^3 (-2 + A)^6 (-1 + A)^{10} \cdot\\
	&A^{14} (1 + A)^{17} (2 + A)^{17} (3 + A)^{14} (4 + A)^8 (16A-80+9B).
\end{align*}

Therefore $\det M=0$ for all $k>0$ is equivalent to $(12 A-48+7B)=0$ or $ (16A-80+9B)=0$ respectively. In either case, it is equivalent to $d'\cdot s_2\in \Z$ by direct calculation.

\section{Proof of Proposition \ref{period}}\label{2}
For every real number $x$, $\{x\}:=x-\fl{x}$ is the fraction part of $x$. 
\begin{lemma}\label{sd}
	For any two real numbers $x,y$, integer $n$, we have:
	\begin{enumerate}
	\item 	$\fl{x}+\fl{y}\leq\fl{x+y}\leq \fl{x}+\fl{y}+1$;
	\item 	$\cl{x}+\cl{y}-1\leq\cl{x+y}\leq \cl{x}+\cl{y}$;
	\item	$n\leq x \ifff n\leq \fl{x}$;
	\item	$n>x \ifff n\geq \fl{x}+1$.
	\end{enumerate}
\end{lemma}
\pf.  Obvious from definition.\qed

\begin{lemma}
\label{lmphi}
	Let $x,y$ be real numbers. Define $\phi(x,y):=\fl{x+y}-\cl{x}+1$, then
	$\fl{y}\leq \phi(x,y)\leq \fl{y}+1.$
\end{lemma}
\pf.
When $x\in\Z$, by Lemma \ref{sd}, $\phi(x,y)=x+\fl{y}-x+1=\fl{y}+1$.
When $x\not\in\Z$, Lemma \ref{sd} shows that
$\fl{x+y}=\fl{x}+\fl{y}$ or $\fl{x}+\fl{y}+1$. Now $\cl{x}-\fl{x}=1$, hence $\phi(x,y)=\fl{y}$ or $\fl{y}+1$. \qed 

\begin{lemma}	\label{clm}
	Consider rational numbers $s_1<s_2<s_3$ with width $w<1$.
	Let $l_k=\fl{k s_2}-\cl{k s_1}+1$, and $r_k=\fl{k s_3}-\cl{k s_2}+1$. Then for every integer $k\geq 1$ we have:
	\begin{align*}
	\fl{k(s_2-s_1)}&\leq l_k\leq 	\fl{k(s_2-s_1)}+1,\\
	\fl{k(s_3-s_2)}&\leq r_k\leq 	\fl{k(s_3-s_2)}+1.
	\end{align*}
	Furthermore, both sequences $\{l_k\}_k$ and $\{r_k\}_k$ are increasing and positive. 
\end{lemma}
\pf.  Let $\phi$ be as defined in Lemma \ref{lmphi}. Then $l_k=\phi(ks_1,k(s_2-s_1))$, and $r_k=\phi(ks_2,k(s_3-s_2))$. Hence the inequalities follows.

Next, since $w<1$, we have $s_2-s_1>1$. Therefore $l_{k+1}\geq \fl{(k+1)(s_2-s_1)}\geq \fl{k(s_2-s_1)}+\fl{s_2-s_1}\geq \fl{k(s_2-s_1)}+1\geq l_k$. So $\{l_k\}_k$ is increasing. Similarly, since $s_3-s_2>1$, the sequence $\{r_k\}_k$ is increasing. In particular, let $k=1$. Then $l_1\geq \fl{s_2-s_1}\geq 1$ and $r_1\geq \fl{s_3-s_2}\geq 1$.\qed

\begin{prop}
	Consider rational numbers $s_1<s_2<s_3$ with width $w<1$. Then for all $n\geq 0$, we have $\pi(n)<w(n+1)$. \label{pinw}
\end{prop}
\pf. We define two auxiliary functions $\pi^-(n)$ (and $\pi^+(n)$) as the number of positive integer $k$ such that $l_k \leq n$ ($r_k\leq n$ respectively). Then by definition, $\pi^-(n)+\pi^+(n)=\pi(n)$.

Next we find bounds for $\pi^-(n)$. Since $\{l_k\}_k$ is increasing by Lemma \ref{clm}, $\pi^-(n)$ equal to the unique integer $k$ such that $l_k\leq n<l_{k+1}$, or equivalently,	$l_k\leq n\leq l_{k+1}-1$.
By Lemma \ref{clm}, we have
	$\fl{k(s_2-s_1)}\leq n\leq (\fl{(k+1)(s_2-s_1)}+1)-1$, i.e.,
$\fl{k(s_2-s_1)}\leq n\leq \fl{(k+1)(s_2-s_1)}.$
Applying Lemma \ref{sd} we find
\begin{align*}
	n&\geq \fl{k(s_2-s_1)} \ifff n+1>k(s_2-s_1),\\
	n&\leq \fl{(k+1)(s_2-s_1)} \ifff n\leq (k+1)(s_2-s_1).
\end{align*}
Hence $k(s_2-s_1)-1 < n \leq (k+1)(s_2-s_1)$. Now  we can solve for $k=\pi^-(n)$:
\begin{equation}\frac{n}{s_2-s_1}-1\leq \pi^-(n)<\frac{n+1}{s_2-s_1}.	\label{npi}\end{equation}
Similar argument for $\{r_k\}$ and $\pi^+(n)$ shows that
\begin{equation}	\frac{n}{s_3-s_2}-1\leq \pi^+(n)<\frac{n+1}{s_3-s_2}.	\label{mpi}\end{equation}
Adding (\ref{npi}) and  (\ref{mpi}), and noticing that $ w=(s_2-s_1)\ik+(s_3-s_2)\ik$, we find
\[wn-2\leq \pi(n)< w(n+1).\label{keypi}\]\qed

\pf {\;\em of Proposition \ref{period}.}  Part (1) follows from Proposition \ref{pinw}. Indeed $\pi(n)\leq w(n+1)<n+1$. 
To prove (2), it suffices to show
\begin{equation}
	\pi(n+m)=\pi(n)+mw.
	\label{nmw}
\end{equation}
for every $n\geq 1$. 

Indeed, assume this holds, then by $w<1$, and part one, we must have $\pi(n+m)< \pi(n)+m\leq n+m$. Therefore, if $n>m$, then $\pi(n)<n$.
It remains to prove for $mw\leq n\leq m$. First let $n=1$ in (\ref{nmw}), then $\pi(m+1)=\pi(1)+mw$.
By the first part of Theorem \ref{period}, $\pi(1)=0$ or $1$. 
Suppose there exists some $n'$ such that $mw\leq n'\leq m$ and  $\pi(n')=n'$. By the definition of $\pi(n)$, $\pi$ in non-decreasing, so
\[\pi(n')\leq \pi(m+1).\]
Hence there are two possible cases:

1). If $\pi(1)=0$, then $n'=\pi(n')\leq \pi(m+1)=mw$. Because $mw\leq n'$, we find $n'=mw$. Hence $\pi(mw)=\pi(m+1)=mw$. However, the column of $\Delta_1$ on $x=0$ contains exactly $m+1$ lattice points. Therefore, there exists $i,j\geq 1$ such that 
$l_i=r_j=m+1$. Hence for any $z<m+1$, $\pi(z)\leq \pi(m+1)-2$. In particular, $\pi(mw)\leq \pi(m+1)-2=mw-2$, a contradiction.

2). If $\pi(1)=1$, then $n'=\pi(n')\leq \pi(m+1)=mw+1$. Because $mw\leq n'$, we find $n'=mw$ or $mw+1$. If $n'=mw$, then $\pi(mw)=mw$. If $n'=mw$, then $\pi(mw+1)=mw+1$. Now the same argument in 1) shows $\pi(mw)\leq \pi(m+1)-2$ and $\pi(mw+1)\leq \pi(m+1)-2$, where $\pi(m+1)-2=mw-1$. Hence we reach a contradiction.

Therefore, no such $n'$ exists. By part (1) of the Theorem again, for every $mw\leq n\leq m$, $\pi(n)<n$.

Now we prove (\ref{nmw}) in the following. By definition, $m$ is the smallest integer such that $m\Delta_0$ is a good lattice triangle. Recall that the coordinates $p=(x_1,y_1)$, $q=(x_2,y_2)$, and the $y$-intercept is $y_0=1$. Hence $m$ equals to the lowest common denominator of $\{x_1,x_2,y_1,y_2\}$.
Let $s_k=u_k/v_k$ be in the lowest terms of $s_k$ for $i=1,2,3$. Let $\alpha=\lcm(v_1,v_2)$, and $\beta=\lcm(v_2,v_3)$. Then $\alpha s_1$ and $\alpha s_2$ are both integers, and $\gcd(\alpha s_1, \alpha s_2, \alpha)=1$ by the choice of $\alpha$. 
We let $m_1$ be the lowest common denominator of $x_1$ and $y_1$, then $m_1\mid \alpha s_2-\alpha s_1$, so there exists an integer $t_1$ such that $m_1 t_1=\alpha s_2-\alpha s_1$.
We claim $t_1\mid \alpha$ and $t_1\mid \alpha s_1$.
Indeed, let $g=\gcd(\alpha s_2-\alpha s_1, \alpha)$, then $(\alpha s_2-\alpha s_1)/g$ equals to the denominator of the lowest term of $x_1$, hence it divides $m_1$. Since $(\alpha s_2-\alpha s_1)/t_1=m_1$, we have $t_1\mid g$, and hence $t_1\mid \alpha$. Similarly, $t_1\mid \alpha s_1$.
Note that $t_1\mid \alpha s_2-\alpha s_1$, hence $t_1\mid \alpha s_2$. However, $\gcd(\alpha s_1, \alpha s_2, \alpha)=1$, so $t=1$, and $m_1=\alpha s_2-\alpha s_1$. 

Similarly, the lowest common denominator $m_2$ of $x_2$ and $y_2$ equals to $\beta s_2-\beta s_3$.
Hence $m=\lcm (m_1,m_2)=\lcm( \alpha (s_2-s_1), \beta(s_3-s_2))$. 

Now for every $k>0$, $l_{k+\alpha}=\fl{(k+\alpha)s_2}-\cl{(k+\alpha)s_1}+1=\fl{ks_2}-\cl{ks_1}+1 +\alpha s_2-\alpha s_1=l_k+\alpha(s_2-s_1)$.
Because the sequence $\{l_k\}$ is non-decreasing, for any given positive integer $n\geq 1$, there exists a unique $k$ such that $l_k\leq n<l_{k+1}$. By the definition of $\pi^+(n)$ (see proof of Proposition \ref{pinw}), $\pi^+(n)=k$. Adding  $\alpha(s_2-s_1)$ to the inequality, we find
\[
	l_{k+\alpha}=l_k+\alpha(s_2-s_1)\leq n+\alpha(s_2-s_1) < l_{k+1}+\alpha(s_2-s_1)=l_{k+1+\alpha}.
\]
This implies that $\pi^+ (n+\alpha(s_2-s_1))=k+\alpha=\pi^+(n)+\alpha$ for all $n$.

Similarly, $r_{k+\beta}=r_k+\beta(s_3-s_2)$ for every $k> 0$, and  $\pi^-(n+\beta(s_3-s_2))=\pi^-(n)+\beta$ for all $n>0$.

Finally, recall $m=\lcm( \alpha (s_2-s_1), \beta(s_3-s_2))$. If we let $t=\gcd(\alpha (s_2-s_1), \beta(s_3-s_2))$, then $mt=\alpha (s_2-s_1)\beta(s_3-s_2)$.
By iteration, we find
\[	\pi^+(n+m)=\pi^+(n+\alpha (s_2-s_1)\cdot \frac{\beta(s_3-s_2)}{t})=\pi^+(n)+\frac{\alpha\beta(s_3-s_2)}{t}.\]
and 
\[\pi^-(n+m)=\pi^-(n+\beta(s_3-s_2)\cdot \frac{\alpha (s_2-s_1)}{t}))=\pi^-(n)+\frac{\alpha\beta(s_2-s_1)}{t}.\]
Adding them up, we conclude
\[
	\pi(n+m)=\pi(n)+\frac{\alpha\beta (s_3-s_1)}{t}=\pi(n)+mw.
\] \qed

\section{Proof of Proposition \ref{main02}}\label{4}
We first point out that for rational numbers $s_1<s_2<s_3$ with width $w<1$, the reduced degree $d=0$ if and only if the minimal degree $d'=0$. Let $ l_k=\fl{k s_2}-\cl{k s_1}+1$, and $ r_k=\fl{k s_3}-\cl{k s_2}+1$ as in Theorem \ref{main03}. Let $\pi(n)$ be defined as in \ref{pi}. Then the following are all equivalent. 
\begin{enumerate}
	\item $d'=1$;
	\item $\pi(1)=1$;
	\item Exactly one of $l_1$ and $r_1$ equals one. 
\end{enumerate}
Clearly, $(3)\implies (1) \iff (2)$ by Proposition \ref{period}. By Lemma \ref{clm}, $l_1\geq 1$ and $r_1\geq 1$. Therefore if $\pi(1)=1$, then one of $l_1$ and $r_1$ is one but not both. So $(2)\implies (3)$.

\begin{lemma}
	Given rational numbers $s_1<s_2<s_3$ with width $w<1$. For all positive integers $k$, define $ l_k=\fl{k s_2}-\cl{k s_1}+1$, and $ r_k=\fl{k s_3}-\cl{k s_2}+1$. Then:
 \begin{enumerate}
		\item If $l_1\neq 1$, then the sequence $\{l_k\}$ is strictly increasing.
		\item If $r_1\neq 1$, then the sequence $\{r_k\}$ is strictly increasing. 
	\end{enumerate}
	\label{strictI}
\end{lemma}

\pf. We need only prove this for $\{r_k\}$. Since $w<1$, $r_{k+1}\geq r_k$ for any $k$. Suppose there exists $k\geq 1$ such that $r_k=r_{k+1}=d+1$. Then by definition we find $
	\fl{ks_3}-\cl{ks_2}=d$, and 
	$\fl{(k+1)s_3}-\cl{(k+1)s_2}=d.$
By adding the same integer to $s_i$, for $i=1,2,3$, we can assume $0< s_2\leq 1$, without changing the values of $\{r_k\}$. Hence $\cl{s_2}=1$. Further,
$r_1\neq 1$, and $w<1$ implies that $r_1\geq 1$, so $r_1\geq 2$. This implies $s_3\geq 2$. By Lemma \ref{clm} we have
	$\cl{(k+1)s_2}\leq \cl{ks_2}+\cl{s_2}=\cl{ks_2}+1$.
Hence	
\[
\fl{(k+1)s_3}=\fl{(k+1)s_2}+d\leq \cl{ks_2}+1+d=\fl{ks_3}+ 1.
\]
On the other hand, $\fl{(k+1)s_3}\geq \fl{ks_3}+\fl{s_3}\geq \fl{ks_3}+2$. So we reached a contradiction. Hence $\{r_k\}$ is strictly increasing. \qed

Let $\pi(n)$ be as defined in Section \ref{mainresult}. Let $\delta(n):=\pi(n)-n$. By Proposition \ref{period}, $\delta(n)\leq 0$ for all $n\geq 1$.
\begin{lemma}\label{strictII}
	Given rational numbers $s_1<s_2<s_3$ with width $w<1$.
	Suppose $\pi(1)=0$ (so that $\delta(1)=-1$). Then for any $n>1$, the following two are equivalent:
	\begin{enumerate}
		\item $\pi(n)=n$ (so that $\delta(n)=0$).
		\item There exists a positive integer $n_0\leq n$ such that $\delta(n_0)=\delta(n_0+1)=\dots=\delta(n)=0$, and $n_0=r_v=l_u$ for some $u,v\geq 1$.
	\end{enumerate}
\end{lemma}
\pf. $\pi(n)=n$ if and only if $\delta(n)=0$. Since $\pi(1)=0$, $l_1\neq 1$ and $r_1\neq 1$. By Proposition \ref{strictI}, both $\{r_j\}$ and $\{l_i\}$ is strictly increasing. Hence for any given integer $m$, there are at most one $r_j$ and at most one $l_i$ equaling to $m$.

The sufficiency in the corollary is clear from definition. So we prove the necessity.
Suppose for a given $n>1$, $\delta(n)=0$, then there are four cases:
\begin{itemize}
	\item There exists $v,u\geq 1$ such that $n=r_v=l_u$. Then $\pi(n-1)=\pi(n)-2=n-2$. Hence $\delta(n-1)=\pi(n-1)-(n-1)=-1$. This shows $n_0=n$.
	\item There exists $v\geq 1$ such that $n=r_v$, and no $l_u=n$. Then $\pi(n-1)=\pi(n)-1=n-1$. Hence $\delta(n-1)=0$. This shows we can reduce the argument from $n$ to $n-1$.
	\item  There exists $l\geq 1$ such that $n=l_u$, and no $r_v=n$. Then $\pi(n-1)=\pi(n)-1=n-1$. Hence $\delta(n-1)=0$. This shows we can reduce the argument from $n$ to $n-1$.
	\item $n$ does not equals to any elements in $\{l_v\}$ and $\{r_v\}$. In this case $\pi(n-1)=\pi(n)=n$, which contradicts to $\pi(n-1)\leq n-1$, hence impossible.
\end{itemize}

As a result, when $n_0\neq n$, we can run the argument again for $n-1$. We claim that eventually it will terminate at some $n_0\geq 2$. Indeed, if it stops at $1$, then there exist $v,u\geq 1$ such that $r_v=l_u=1$, which contradicts to the assumption that $\pi(1)=0$. If it reduces to $0$, then there exist $v\geq 1$ or $u\geq 1$ such that $r_v=0$ or $l_u=0$, which contradicts to $w<1$. As a result, $n_0=l_{u'}=r_{v'}$ for some positive integers $u',v'$ and $n_0\geq 2$.
This finishes the proof of necessity.\qed

\begin{corollary}	\label{crossref}
	Given rational numbers $s_1<s_2<s_3$ with width $w<1$. Assume $\pi(1)=0$. Then the minimal degree $d'\geq 2$ if and only if there exists a positive integer $v$ such that 
	\[r_v=l_{r_v-v}.\]
\end{corollary}
\pf. Since $\pi(1)=0$, by Lemma \ref{strictI} and Lemma \ref{clm}, both $\{l_k\}$ and $\{r_k\}$ are strictly increasing and at least $2$. If $d'\geq 2$, then there exists $n\geq 2$ such that $\pi(n)=n$. Then by Lemma \ref{strictII}, there exists $n_0\leq n$ and $u,v\geq 1$, such that $n_0=r_v=l_u$ and $\pi(n_0)=n_0$. Since both $\{l_k\}$ and $\{r_k\}$ are strictly increasing, so $\pi(n_0)=u+v$ by definition of $\pi$. Now $u+v=\pi(n_0)=n_0=r_v$, so $u=r_v-v$.

Conversely, suppose $r_v=l_{r_v-v}$ for some $v\geq 1$. Let $n_0=r_v$. Then $n_0\geq r_1\geq 2$. Since both $\{l_k\}$ and $\{r_k\}$ are strictly increasing, $\pi(n_0)=v+r_v-v=r_v=n_0$. Therefore $d'\geq 2$ by definition of the minimal degree.
\qed

\para

\pf {\;\em of Proposition \ref{main02}}. Suppose $d\geq 1$, we prove that one of  (1), (2), (3) is false. Since $w<1$, both $l_1$ and $r_1$ are positive. Since $d\geq 1$, $d'\geq 1$ too. If $d'=1$, then $l_1=1$ or $r_1=1$, so (1) fails. Otherwise, assume $d'\geq 2$. Now $\pi(1)=0$. By Corollary \ref{crossref}, we conclude that there exists a positive integer $v$ such that  $r_v=l_{r_v-v}$. Let $u(v)=r_v-v$.

Recall $\gamma$ is the smallest positive integer such that $\gamma s_2^2\in\Z$, $\gamma s_3\in\Z$, and $\gamma s_2 s_3\in\Z$. Therefore, $\gamma s_2\in \Z$ also holds. Let  $q$ be the quotient of $v$ by $\gamma$, and $t$ be the remainder. That is, $v=q\gamma +t$ with $0\leq t<\gamma$. Then $q$ and $t$ are unique. We define $r_0=1$.
Now $r_v=\fl{vs_3}-\cl{vs_2}+1=\fl{(q\gamma+t)s_3}-\cl{(q\gamma+t)s_2}+1$, which equals to $\fl{ts_3}-\cl{ts_2}+1+q\gamma(s_3-s_2)=r_t+q\gamma (s_3-s_2)$ as $\gamma s_2\in\Z$ and $\gamma s_3\in\Z$. Hence 
\[r_v-v=r_t-t+q\gamma (s_3-s_2-1).\]
Let $A=\gamma (s_3-s_2-1)$, so that $r_v-v=r_t-t+qA$. Then $A>0$ since $w<1$.
Next we calculate $l_u=l_{r_v-v}$. We have
\begin{align*}
	l_u&= \fl{us_2}-\cl{us_1}+1\\
	&=\fl{(r_t-t+qA)s_2}-\cl{(r_t-t+qA)s_1}+1\\
	&=\fl{(r_t-t)s_2}+qAs_2-\cl{(r_t-t+qA)(s_1-s_2)+(r_t-t+qA)s_2}+1\\
	&=\fl{(r_t-t)s_2}+qAs_2-\cl{(r_t-t+qA)(s_1-s_2)+(r_t-t)s_2}-qAs_2+1\\
	&=\fl{(r_t-t)s_2}-\cl{(r_t-t+qA)(s_1-s_2)+(r_t-t)s_2}+1.
\end{align*}
where we applied the fact that $qAs_2=q\gamma (s_3 s_2-s_2^2-s_2)$ is an integer.

Replacing $l_u=r_v$ by the expressions of them above, we find
\begin{align*}
	\fl{(r_t-t)s_2}-\cl{(r_t-t+qA)(s_1-s_2)+(r_t-t)s_2}+1&=r_t+q\gamma(s_3-s_2)\\
	\fl{(r_t-t+qA)(s_2-s_1)-(r_t-t)s_2}&=r_t-1+q\gamma(s_3-s_2)-\fl{(r_t-t)s_2}.
\end{align*}
Hence by definition of the floor function,
\begin{align*}
	r_t-1+q\gamma(s_3-s_2)-\fl{(r_t-t)s_2}&\leq  (r_t-t+qA)(s_2-s_1)-(r_t-t)s_2\\
	&<r_t+q\gamma(s_3-s_2)-\fl{(r_t-t)s_2}.
\end{align*}
Adding $(r_t-t)s_2$ to all sides, we obtain that
\[r_t-1+q\gamma(s_3-s_2)+\{(r_t-t)s_2\}\leq  (r_t-t+qA)(s_2-s_1)<r_t +q\gamma(s_3-s_2)+\{(r_t-t)s_2\}.\]
Recall $r_t-t+qA=r_v-v$. Since $\{r_k\}$ is strictly increasing and $r_1\geq 2$, we have $r_k>k$ for all $k\geq 1$. Hence $r_v>v$, so we can divide by $r_t-t+qA$:
\[\frac{r_t-1+q\gamma(s_3-s_2)+\{(r_t-t)s_2\}}{ r_t-t+qA}\leq s_2-s_1<\frac{r_t +q \gamma(s_3-s_2)+\{(r_t-t)s_2\}}{r_t-t+qA}.\]

Let $B=r_t+\{(r_t-t)s_2\}$ and $C=r_t-t$. Then the above is equivalent to
\begin{align}
\label{main}\frac{q(A+\gamma)+B-1}{qA+C}\leq s_2-s_1<\frac{q(A+\gamma)+B}{qA+C}.
\end{align}

Notice that all steps above are equivalent. Consequently, when $d'\neq 1$, for any $v\geq 1$, $l_{r_v-v}=r_v$ if and only if $v=q\gamma+t$ and $\gamma, t$ satisfy (\ref{main}).

\para
Now suppose $d'\geq 2$ and  $l_{r_v-v}=r_v$ for some $v\geq 1$. Then $v=q\gamma+t$ and $\gamma, t$ satisfy (\ref{main}).
Since $w<1$, we have $s_2-s_1>\dis\frac{s_3-s_2}{s_3-s_2-1}=(A+\gamma)/A$. As a result,
\[\frac{A+\gamma}{A}<\frac{q(A+\gamma)+B}{qA+C}.\]
Since $A,B,C>0$, we find $(A+\gamma)C<AB$. 
Therefore, if $t\neq 0$, then $q\geq 0$, so 
\[	s_2-s_1<\frac{q(A+\gamma)+B}{qA+C}\leq \frac{B}{C}.\]
which shows that (2) of Proposition \ref{main02} is false. If $t=0$, then $q>0$ because $v=q\gamma+0>0$, so
\[	s_2-s_1<\frac{q(A+\gamma)+B}{qA+C}\leq \frac{(A+\gamma)+B}{A+C}.\]
Here $t=0$, so $r_0=1$, $B=1+\{s_2\}$ and $C=1$. Therefore
\[s_2-s_1<\frac{\gamma(s_3-s_2)+1+\{s_2\}}{\gamma (s_3-s_2-1)+1}.\]
which says (3) of Proposition \ref{main02} is false.
This proves the first part of the proposition.

\para
Conversely, assume $(4)$ of Proposition \ref{main02}  holds. If (1) fails then $d'=1$, so $d\neq 0$. Otherwise, suppose (1) holds (so $d'\neq 1$) and one of (2) and (3) is false, we will prove that $d'\geq 2$, so $d\neq 0$.

{\bf Case I.}
Assume (1) holds and (2) is false. Let $t$ be the integer such that $1\leq t\leq \gamma-1$ and
\[\frac{r_t+\{(r_t-t)s_2\}}{r_t-t}=\max_{1\leq i\leq \gamma-1} \frac{r_i+\{(r_i-i)s_2\}}{r_i-i},\]
Then by $w<1$ we have
\[	\frac{s_3-s_2}{s_3-s_2-1}<s_2-s_1<\frac{r_t+\{(r_t-t)s_2\}}{r_t-t},\]
which is
\[\frac{A+\gamma}{A}<s_2-s_1<\frac{B}{C}.\]
So  $(A+\gamma)C<AB$.
Define two sequences $\{x_n\}$ and $\{y_n\}$ by
\[x_n=\frac{n(A+\gamma)+B-1}{nA+C}, \quad y_n=\frac{n(A+\gamma)+B}{nA+C}.\]
Then $x_n< y_n$ for all $n\geq 0$. Further
\[	\lim_{n\ra \infty} x_n=\lim_{n\ra \infty} y_n=\frac{s_3-s_2}{s_3-s_2-1}=\frac{A+\gamma}{A}.\]
Since $y_{n+1}-y_n=(nA+C)\ik ( (n+1)A+C)\ik ((A+\gamma)C-AB)$, and $(A+\gamma)C<AB$, $\{y_n\}$ is strictly decreasing,
We now prove the following interval inclusion:
\begin{align}\label{inclusion}
\big(\frac{A+\gamma}{A},\frac{B}{C}\big)\subset \bigcup_{n\geq 0} [x_n, y_n).
\end{align}
Then there exists a non-negative integer $q$ such that $x_q\leq s_2-s_1< y_q$, which is exactly the equation (\ref{main}). Let $v=q\gamma+t$, then $v\geq 1$ and $l_{r_v-v}=r_v$.  Hence $d'\geq 2$ by Corollary \ref{crossref}. 

So we prove the inclusion (\ref{inclusion}) above. Indeed, $\{y_n\}$ is strictly decreasing, so we need only to show that $x_n\leq y_{n+1}$ for all $n\geq 0$. We have
\[x_q-y_{q+1}=\frac{n(A+\gamma)+B-1}{nA+C}-\frac{(n+1)(A+\gamma)+B}{(n+1)A+C}.\]

Since $w<1$, $A,C>0$. Hence $x_n-y_{n+1}\leq 0$ if and only if
\[(n(A+\gamma)+B-1)((n+1)A+C)-( (n+1)(A+\gamma)+B)(nA+C)\leq 0.\]
which simplifies to
\[A(B-n-1)\leq (A+\gamma+1)C.\]
Now (4) of  Proposition \ref{main02} says $A(B-1)\leq (A+\gamma+1)C$. Since $A>0$, $n\geq 0$, $A(B-n-1)\leq (A+\gamma+1)C$ follows.

{\bf Case II.} Suppose (1) holds and (3) is false. Then by $w<1$ we have
\[\dis\frac{s_3-s_2}{s_3-s_2-1}<s_2-s_1<\frac{\gamma (s_3-s_2)+1+\{s_2\}}{\gamma(s_3-s_2-1)+1}.\]
which is
\begin{align}\label{t0case}
\frac{A+\gamma}{A}<s_2-s_1<\frac{A+\gamma+1+\{s_2\}}{A+1}.
\end{align}
So  $(A+\gamma)<A(1+\{s_2\})$, i.e., $\gamma <A\{s_2\}$.

Consider two sequences $z_n$ and $w_n$ defined by
\[z_n=\frac{n(A+\gamma)+\{s_2\}}{nA+1}, \quad w_n=\frac{n(A+\gamma)+1+\{s_2\}}{nA+1}.\]

Since $\{s_2\}<1<(s_3-s_2)/(s_3-s_2-1)$, we have
$z_n<(A+\gamma)/A<w_n$ for all $n\geq 0$. 
Furthermore, $w_{n+1}-w_n=(nA+1)\ik ((n+1)A+1)\ik (\gamma-A\{s_2\})$. So $\gamma<A\{s_2\}$ implies that $\{w_n\}$ is strictly decreasing.
Finally, $\lim_{n\ra\infty}w_n=(A+\gamma)/A$. 
Now (\ref{t0case}) says $\lim_{n\ra\infty}w_n<s_2-s_1<w_1$. Hence, there exists a positive integer $q$ such that $z_q<(A+\gamma)/A<s_2-s_1<w_q$.

Note that for $t=0$, $r_t=1$, so $B=1+\{s_2\}$ and $C=1$. Hence the equation $z_q<s_2-s_1<w_q$ is the equation (\ref{main}) with $t=0$. Therefore $v=q \gamma+0$ is a positive integer and $l_{r_v-v}=r_v$. Hence $d'\geq 2$.
\qed

\section{Proof of Theorems \ref{main03} and \ref{main04}}
\label{03pf}

We have shown that the first case of $d'=5$ in Theorem \ref{main03} gives a blow-up which is not a MDS if $5s_2\not\in\Z$. The result for $d'=7, 9$ are given by using a computer program. Now to finish the proof of Theorem \ref{main04}, we need only to show that if $2\leq d'\leq 9$ and $d'\cdot s_2\not\in\Z$, then
\begin{enumerate}
	\item Either $\Delta_1$ belongs to Gonz\'{a}lez and Karu's nonexamples (see Remark \ref{gknonex}), or
	\item $d'=5,7$ or $9$, and the slopes satisfy the inequalities in Theorem \ref{main03}.
\end{enumerate}
Then Theorem \ref{main03} and Theorem \ref{main04} will hold.

Consider the triangle $\Delta_1$ of given slopes $(s_1,s_2,s_3)$. We assume $d'\geq 2$.
Since $d'\neq 1$, Lemma \ref{strictI} shows that the numbers of points on each columns are strictly increasing from the leftmost and rightmost vertices to the center column. Notice that there are exactly $d'$ columns in $\Delta_1$ with $\leq d'$ points. Hence these columns must have $2,3,\cdots,d'-1,d',d'$ points respectively. 

Now it is clear that every such triangle $\Delta_1$ determines a subset $S=\{u_1,\cdots,u_{s-1}\}\subset \{2,\cdots, d'-1\}$, whose complement is $\{v_1,\cdots, v_{t-1}\}$, such that leftmost columns of $\Delta_1$ have $u_1,\cdots,u_{s-1},d'$ points, and the rightmost columns have $v_1,\cdots,v_{t-1},d'$ points. In particular, Gonz\'{a}lez and Karu's nonexamples correspond to $S=\emptyset$ and $T=\{2,\cdots,d'-1\}$.

Conversely, every such subset $S$ determines a system of inequalities on the three slopes, by requiring that the leftmost $s$ columns and rightmost $t$ columns have the given number of points, from $S$, $T$, and $d'$. Therefore, to find all possible triangles with the given $d'$, it suffices to solve the systems of inequalities for all subsets of $\{2,\cdots,d'-1\}$ and collect those who have solutions. This job is best done by a computer program. Here we will show an alternative proof for $2\leq d'\leq 6$.

\begin{lemma}
	Let $\Delta$ be a lattice triangle given by slopes $s_1<s_2<s_3$. Let $l_k$ be the number of lattice points on the $k$-th column from the left. 
	Let the left vertex be $p$. Let $n\in \Z_{>0}$. Suppose there are at least $nk$ columns with $x$-coordinate $<0$. Then
	\begin{enumerate}
		\item $l_k-l_{k-1}\geq l_1-1$;
		\item $l_k\geq kl_1-(k-1)$.
		\item $l_{nk}\geq kl_n-(k-1)$.
	\end{enumerate}
Similarly, let $r_k$ be the number of lattice points on the $k$-th column from the right. Suppose there are at least $nk$ columns with $x$-coordinate $>0$. Then the same inequalities hold for $r_k$.
	\label{kstep}
\end{lemma}
\pf. By the hypothesis, $l_k$ has the following expression:
\[l_k=\fl{k s_2}-\cl{k s_1}+1.\]
Hence 
\[l_{k}-l_{k-1}=\fl{k s_2}-\fl{(k-1)s_2}-(\cl{k s_1}-\cl{(k-1)s_1}).\]
By Lemma \ref{sd},
\[\fl{k s_2}-\fl{(k-1)s_2}\geq \fl{s_2}, \quad 	\cl{k s_1}-\cl{(k-1) s_1}\leq \cl{s_1} .\]
Therefore $l_k-l_{k-1}\geq \fl{s_2}-\cl{s_1}=l_1-1$. Now adding the inequalities in $(1)$ gives $(2)$. Finally, notice that the proof of $(1)$ works when replacing the sequence $\{l_k\}_k$ by $\{l_{nk}\}_k$. Hence, $(3)$ holds. The proof for $\{r_k\}$ is identical.\qed

\para
We can now exclude many subsets $S$ and $T$ for which the conditions in Lemma \ref{kstep} cannot be satisfied. As in the lemma, we let $l_k$ ($r_k)$ be the number of lattice points on the $k$-th column from the left (right).

Now we prove the classification part of Theorem \ref{main04} for $d'\leq 6$. Firstly by symmetry we can always assume $\vv{S}\leq \fl{d'/2-1}$.
For $d'=2,3$, the only possible case is $S=\emptyset$.
For $d'=4$, $S=\emptyset$ or $\{2\}$ or $\{3\}$. If $S=\{2\}$, then $T=\{3\}$, so that the $r_1=3$, $r_2=4$. Now $r_{2}<2r_1-1$ contradicts to Lemma \ref{kstep}.
If $S=\{3\}$, then $T=\{2\}$, so by symmetry this case does not exist.

For $d'=5$, if $S=\{2\}$ then $T=\{3,4\}$. Then $(r_1=3,r_2=4,r_3=5)$. If $S=\{4\}$ then $l_1=4,l_2=5$. In both case, $r_{2}<2r_1-1$, contradicts to Lemma \ref{kstep}.
  So $S=\emptyset$ or $S=\{3\}$.

For $d'=6$, Either $2\in S$ or $2\in T$. Assume $2\in S$. Then $l_1=2$, and $r_1\geq 3$. 
If $r_1=4$ or $5$, then by Lemma \ref{kstep}, $r_2\geq 7$, contradiction.
If $r_1=6$, then $T=\emptyset$, contradicts to $\vv{S}\leq 2$.
If $r_1=3$, then either $r_2=5,r_3=6$, or $r_2=6$. The first gives a contradiction since $r_3-r_2\leq 2$. The second gives $S=\{2,4,5\}$, contradicts to $\vv{S}\leq 2$. Therefore $2\not\in S$. 

Now we have $2\in T$. Then $3\leq l_1\leq 6$. If $l_1=6$ then $S=\emptyset$. Otherwise $l_1=3$ or $4$ or $5$. In all cases, $l_3\geq 3l_1-2\geq 7$, so $l_2=6$. But this shows $S=\{3\}$, and $T=\{2,4,5\}$. Then $r_2=4$ and $r_4=6$ contradicts to $(3)$ of Lemma \ref{kstep}.

In conclusion, when $2\leq d'\leq 6$, the only possible triangles are Gonz\'{a}lez and Karu's nonexamples, and the one where $S=\{3\}$, $T=\{2,4\}$. We can solve the corresponding system of inequalities given by $l_1=3,l_2=5$, $r_1=2$, $r_2=4$, and $r_3=5$. Assuming that $0\leq s_2<1$, the solution follows form Lemma \ref{shape}:
\[\begin{cases}
	-2-\frac{1}{2}<s_1\leq -2;\\
	\frac{1}{3}<s_2<\frac{1}{2};\\
	2\leq s_3<2+\frac{1}{3}.
\end{cases}\]

The rest cases for $d'=7,8,9$ are calculated by computer programs. Note that when $d'\leq 6$ we did not invoke the condition of $w<1$. However we need the condition $w<1$ when $d'$ is at least $7$. 
Assume $7\leq d' \leq 9$, and $\vv{S}\leq \fl{d'/2-1}$. The only possible triangles other than the $S=\emptyset$ case are given by

\begin{enumerate}[(a)]
	\item $S=\{4\}$, $T=\{2,3,5,6\}$.
	\item $S=\{3,5\}$, $T=\{2,4,6\}$.
	\item $S=\{3,6\}$, $T=\{2,4,5,7\}$.
	\item $S=\{5\}$, $T=\{2,3,4,6,7,8\}$.
	\item $S=\{3,5,7\}$, $T=\{2,4,6,8\}$.
\end{enumerate}
We show that (a), (c) and (d) do not satisfy the assumptions of Theorem \ref{main03}. By a computer program, the slopes in case (a) satisfies the following inequalities up to adding a same integers to all the slopes:
\[\begin{cases}
		-\frac{5}{2}<s_1\leq -2;\\
	s_2=1;\\
	\frac{7}{3}\leq s_3<\frac{12}{5}.
\end{cases}\]
Therefore the width $w$ is at least the width of triple $(-5/2, 1, 12/5)$, which equals to $1$. This contradicts to our assumption that $w<1$.

Similarly, the slopes in case (d) satisfies
\[\begin{cases}
	-\frac{7}{2}<s_1\leq -3;\\
	s_2=1;\\
	\frac{9}{4}\leq s_3<\frac{16}{7}.
\end{cases}\]
so that the width is at least $1$, a contradiction.

The slopes in case (c) satisfies
\[\begin{cases}
	-\frac{7}{3}<s_1\leq -2;\\
	s_2=\frac{1}{2};\\
	2\leq s_3<\frac{11}{5}.
\end{cases} \txt{or}\quad 
\begin{cases}
	-\frac{5}{3}<s_1\leq -\frac{3}{2};\\
	s_2=1;\\
	\frac{5}{2}\leq s_3<\frac{13}{5}.
\end{cases} 
\]

The second system of inequalities gives $w\geq 1$, hence a contradiction.
The first system of inequalities gives $d'\cdot s_2=4\in\Z$, so it does not satisfy the hypothesis of Conjecture \ref{conj00}.
In conclusion, when assuming $d'\cdot s_2\not\in\Z$, the only non-examples of MDS are given by
\begin{itemize}
	\item $S=\{3\}$, $T=\{2,4\}$, $d'=5$;
	\item $S=\{3,5\}$, $T=\{2,4,6\}$, $d'=7$;
	\item $S=\{3,5,7\}$, $T=\{2,4,6,8\}$, $d'=9$.
\end{itemize}
which proves Theorem \ref{main03} and Theorem \ref{main04}.

\para
\pf {\em\;of Lemma \ref{shape}.} Using the notations $l_k$ and $r_k$ (Lemma \ref{kstep}), the hypothesis says $l_k=2k+1$ for $1\leq k\leq n$. On the right side, $r_k=2k$  for $1\leq k\leq n$ and $r_{n+1}=2n+1$.

We can assume $0\leq s_2<1$. First $l_1=\fl{s_2}-\cl{s_1}+1=3$, and $\fl{s_2}=0$. Hence $\cl{s_1}=-2$. Therefore $\cl{ns_1}\leq -2n$. Now $l_{n}=\fl{ns_2}-\cl{ns_1}+1=2n+1$. So we have $\fl{ns_2}\leq 0$. Since $s_2\geq 0$, the only possibility is $\fl{ns_2}=0$ and $\cl{ns_1}=-2n$. As a result, $ns_2<1$, and $-2n-1<ns_1\leq -2n$. So $s_2<1/n$, and $-2-1/n<s_1\leq -2$.

Next we examine the columns on the right. Since $r_1=\fl{s_3}-\cl{s_2}+1=2$, and $\cl{s_2}=1$, we have $\fl{s_3}=2$, so that $\fl{(n+1)s_3}\geq 2(n+1)$. From $r_{n+1}=\fl{(n+1)s_3}-\cl{(n+1)s_2}+1=2n+1$,  we find $\cl{(n+1)s_2}\geq 2$. Since $s_2<1/n$, the only possibility is $\cl{(n+1)s_2}=2$ and $\fl{(n+1)s_1}=2n+2$. 
Therefore $(n+1)s_2>1$ and $2n+2\leq (n+1)s_1< 2n+3$, so $s_2>1/(n+1)$ and $2\leq s_1< 2+1/(n+1)$.\qed

\appendix 
\section{Figures Classifying triples $(7,b,c)$ where $b,c\leq 70$}

We classify triple $(7,b,c)$ where $7\leq b,c\leq 70$ by various conditions.

\para
In Figure \ref{fig:class1}, we classify the triples by Cutkosky, Gonz\'{a}lez and Karu's results. We first exclude the triples such that $b>c$.
Then we paint the blocks in the following order:
\begin{enumerate}
	\item Those such that $(a+b+c)^2>abc$, so $-K_{X'}$ is big, and hence the blow-up $X'$ is a MDS \cite[Cor. 1]{cutkosky1991}.
	\item Then we exclude triples where $a,b,c$ are not pairwise coprime.
	\item Those in the remaining such that no relation $(e,f,-g)$ with width $w=cg^2/ab<1$ exists, even after permuting $7,b$ and $c$.
	\item Those in the remaining such that $cg<a+b+c$, so $-K_{X'}$ is big (See Remark \ref{big}), and $X'$ is a MDS.
	\item Those in the remaining which belong to Gonz\'{a}lez and Karu's non-examples \cite{GK}. 
	\item The rest.
\end{enumerate}

\para
In Figure \ref{fig:class2}, we classify the triples by their reduced degrees. Firstly we exclude the triples such that $b>c$, or not pairwise coprime. For every triple remaining, if there exists a relation $(e,f,-g)$ such that the width $w=cg^2/ab<1$, then the reduced degree can be defined by Proposition \ref{slopeformula}. Then we divide the triples into three types: $d=0$, $d=1$, or $d\geq 2$. Indeed $d=0$ implies the blow-up is a MDS, and $d=1$ implies the blow-up is not a MDS. For some triples of $d\geq2$ with small minimal degrees, whether the blow-up $X'$ is a MDS is known by Theorems \ref{main03}, \ref{main04}, or by Gonz\'{a}lez and Karu's non-example criterion \cite{GK}.

\newpage
\begin{center}
	\includegraphics[width=1\textwidth]{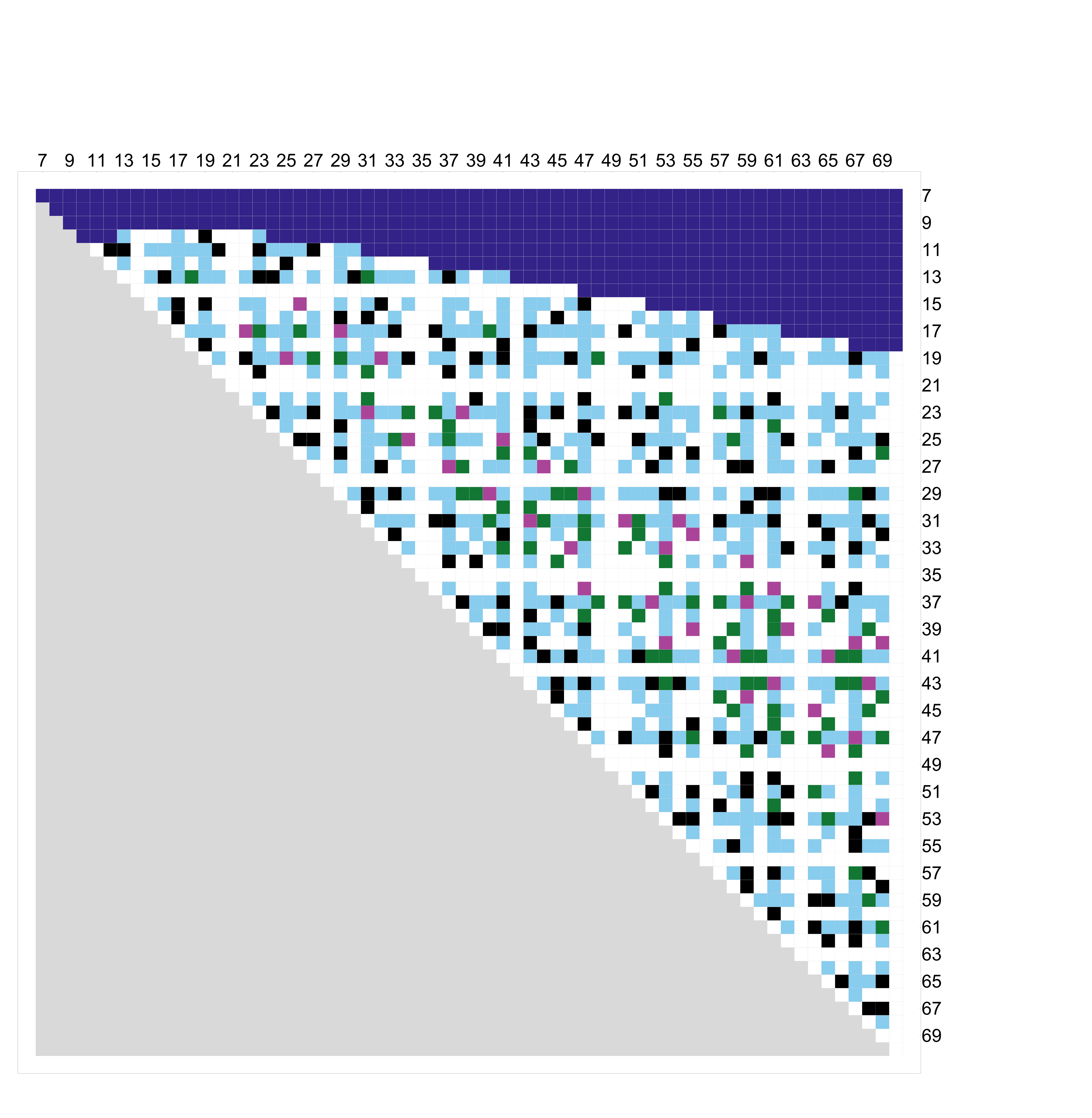}
	\includegraphics[width=0.45\textwidth]{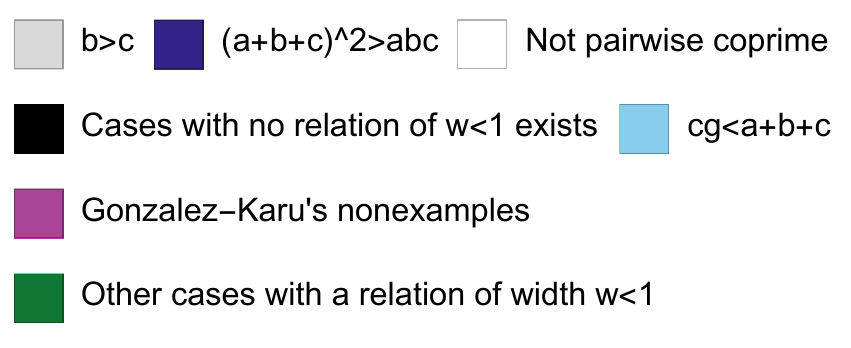}
	\captionof{figure}{Classification of $(7,b,c)$ where $b,c\leq 70$ by Cutkosky, Gonz\'{a}lez and Karu's results.}
	\label{fig:class1}
\end{center}

\begin{center}
	\includegraphics[width=1\textwidth]{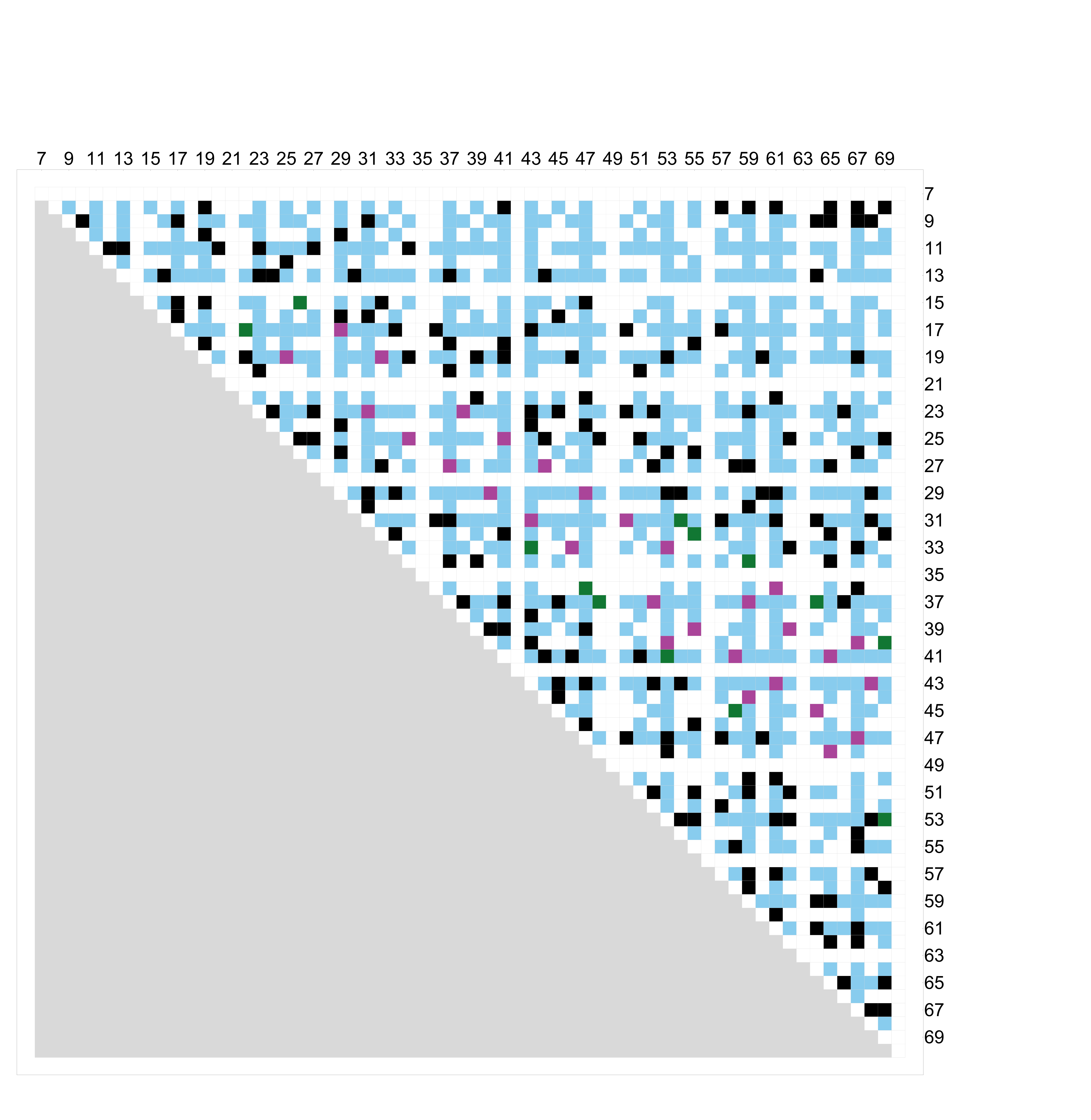}
	\includegraphics[width=0.52\textwidth]{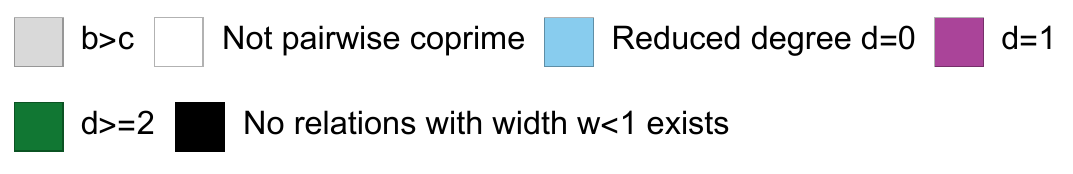}
	\captionof{figure}{Classification of $(7,b,c)$ where $b,c\leq 70$ by reduced degrees}
	\label{fig:class2}
\end{center}

\newpage

\bibliographystyle{halpha}
\bibliography{bibfile}

\begin{thebibliography}{{Wol}16}

\bibitem[Cas16]{castravet16}
Ana-Maria Castravet.
\newblock Mori {D}ream {S}paces and blow-ups.
\newblock {\em Proceedings of the AMS Summer Institute in Algebraic Geometry
  2015}, to appear. Available at
  \href{https://arxiv.org/abs/1701.04738}{arXiv:1701.04738 [math.AG]}, 2016.

\bibitem[CT15]{castravet2015}
Ana-Maria Castravet and Jenia Tevelev.
\newblock $\overline {M}_{0,n}$ is not a {Mori} dream space.
\newblock {\em Duke Mathematical Journal}, 164(8):1641--1667, 06 2015.

\bibitem[Cut91]{cutkosky1991}
Steven~Dale Cutkosky.
\newblock Symbolic algebras of monomial primes.
\newblock {\em J. Reine Angew. Math}, 416:71--89, 1991.

\bibitem[Dum06]{Dumnicki2006}
Marcin Dumnicki.
\newblock {Reduction method for linear systems of plane curves with base fat
  points}, 2006, \href{https://arxiv.org/abs/math/0606716}{arXiv:math/0606716}.

\bibitem[GK16]{GK}
Jos\'{e}~Luis Gonz\'{a}lez and Kalle Karu.
\newblock Some non-finitely generated {C}ox rings.
\newblock {\em Compositio Mathematica}, FirstView:1--13, 2 2016.

\bibitem[GNW94]{goto1994}
Shiro Goto, Koji Nishida, and Kei-ichi Watanabe.
\newblock Non-{C}ohen-{M}acaulay symbolic blow-ups for space monomial curves
  and counterexamples to {C}owsik's question.
\newblock {\em Proceedings of the American Mathematical Society},
  120(2):383--392, 1994.

\bibitem[HK00]{hu2000}
Yi~Hu and Sean Keel.
\newblock Mori dream spaces and {GIT}.
\newblock {\em The Michigan Mathematical Journal}, 48(1):331--348, 2000.

\bibitem[HKL16]{hkl}
J\"uergen Hausen, Simon Keicher, and Antonio Laface.
\newblock On blowing up the weighted projective plane.
\newblock 2016, \href{https://arxiv.org/abs/1608.04542}{arXiv:1608.04542
  [math.AG]}.

\bibitem[KM08]{kollar2008}
Janos Koll{\'a}r and Shigefumi Mori.
\newblock {\em Birational Geometry of Algebraic Varieties}.
\newblock Cambridge University Press, 2008.

\bibitem[Laz04]{lazarsfeldpos}
Lazarsfeld.
\newblock {\em Positivity in Algebraic Geometry I: Classical Setting: Line
  Bundles and Linear Series}.
\newblock Springer, 2004.

\bibitem[Sri91]{HS}
Hema Srinivasan.
\newblock On finite generation of symbolic algebras of monomial primes.
\newblock {\em Communications in Algebra}, 19(9):2557--2564, 1991.

\bibitem[{Wol}16]{mma104}
{Wolfram Research, Inc.}
\newblock Mathematica 10.4.
\newblock \url{https://www.wolfram.com}, Champaign, Illinois, 2016.

\end{thebibliography}

\end{document}